\documentclass[11pt]{article}
\usepackage{graphicx,epsfig,amssymb,cite}
\usepackage{amsmath}

\newtheorem{theorem}{Theorem}
\newtheorem{remark}{Remark}
\newtheorem{lemma}{Lemma}
\newtheorem{example}{Example}
\newtheorem{assumption}{Assumption}
\newtheorem{corollary}{Corollary}

\newtheorem{definition}{Definition}

\date{}

\begin{document}
\title{Mixed $H_-/H_{\infty}$  Fault Detection  Filtering  for It\^o-Type Affine Nonlinear Stochastic Systems}
\author{ Tianliang~Zhang ${}^1$,
 Feiqi~Deng ${}^1$\thanks{
  Email: t\_lzhang@163.com(T. Zhang), aufqdeng@scut.edu.cn(F. Deng), w\_hzhang@163.com(W. Zhang), bschen@ee.nthu.edu.tw(B. S. Chen).
 }, Weihai~Zhang ${}^2$, and Bor-Sen Chen ${}^3$\\
%EndAName
${}^1$ {\small School of Automation Science and Engineering,} \\
 {\small South China University of Technology,
Guangzhou 510640,   P. R. China}\\
${}^2$ {\small College of Information and Electrical Engineering,} \\
{\small Shandong University of Science and Technology,}  {\small
Qingdao  266510, P. R. China}\\
${}^3$ {\small Department of Electrical Engineering,}\\
 {\small National Tsing Hua University, Hsinchu, 30013, Taiwan.}}
\maketitle

{\bf Abstract-} This paper studies the mixed $H_-/H_{\infty}$ fault detection
filtering of It\^o-type nonlinear stochastic systems.  Mixed
$H_-/H_{\infty}$  filtering  combines the system  robustness  to the
external disturbance and the sensitivity to the fault of the
residual signal.  Firstly, for  It\^o-type affine nonlinear
stochastic systems, some sufficient criteria are obtained for the
existence of  $H_-/H_{\infty}$ filter in terms of Hamilton-Jacobi
inequalities (HJIs). Secondly, for  a  class of quasi-linear  It\^o
systems,  a sufficient condition is given for the existence of
$H_-/H_{\infty}$ filter by means of  linear matrix inequalities
(LMIs).  Finally, a numerical example is presented to illustrate the
effectiveness of the proposed results.


\section{Introduction}

Along with the development of modern  industrial production, higher
requirements for safety and reliability  have been put forward. In
order to ensure safety and reliability in industrial process,
various techniques for  fault detection, fault isolation  and fault
estimation have appeared \cite{Stevenbook,Meskinbook}. Fault
detection filter (FDF) is to use the estimated value of the  system state
and the measurement output  to  generate  the residual signal to
detect the system fault. According to the real time comparison
between the designed residual evaluation function and the
corresponding threshold, we are in a position to determine whether
there is a fault occurring. In real world, some unknown disturbance
signals existing in  dynamic systems will fluctuate the residual.
Therefore, one  needs to design robust FDF  to
exclude these unknown external effect. $H_\infty$ control is one of
the most important robust control methods since G. Zames'
fundamental work \cite{G.Zames_81} published,  which  has been
studied extensively
\cite{v.dragan_2005,Petersen,Gershon_05,bshen_13,zhangsiam,zhangfilter,newtsp2,book},
and has been applied to event triggers \cite{wangzhang,newtsp3}, fault
diagnosis \cite{Lih2017,Suxiaojie2016,Zhong2003}, sliding mode control method \cite{newkao1} and adaptive control method \cite{newkao2}. $H_{\infty}$ FDF  is to design the filter such that the $L_2$-gain
from the external  interference to the residual signal is less than
the given attenuation level $\gamma>0$, which reflects the robust
ability of the concerned systems. $H_-$  FDF is to design the filter such that  the $L_2$-gain from the
fault signal to the residual signal is larger than the given
sensitivity level $\delta>0$, which measures  the  sensitivity of
the considered systems to the fault signal. Different from the sole
$H_{\infty}$ or  $H_-$  FDF,  mixed
$H_-/H_{\infty}$ FDF is a combination of
$H_{\infty}$ FDF and $H_-$  FDF, which not only  meets the robustness requirement, but also
requires the residual signal to be sensitive enough to the fault
information
  \cite{Chadli2013,Wangzhenhua2017}. So $H_-/H_{\infty}$ FDF  is one of the most popular robust design methods.

Stochastic  It\^o-type systems are ideal mathematical models in
finance mathematics \cite{Yong_1999}, systems biology
\cite{bschen_2008_1,chench}, benchmark mechanical systems
\cite{xjxie}, so the study for stochastic systems has attracted many
researchers'  interest, and stochastic control has become one of the
most important research fields in modern control theory \cite{newyin1,newyin2,newwu1}.
In 1998,
$H_\infty$ control of linear It\^o systems was  first investigated
in \cite{Hinrichsen,ug}, and from then, nonlinear $H_\infty$ control
\cite{berman_06_1,zhangsiam} and filtering \cite{zhangfilter}, mixed
$H_2/H_\infty$ control \cite{book}  of It\^o systems have   been
solved.  $H_-$   FDF  is to use the $H_-$ index to
measure the minimum influence of the fault on the residual signal
\cite{liuwangyang,liliu,lizhou,Suxiaojie2016,Zhong2003}, where
\cite{liuwangyang,Zhong2003} and \cite{liliu,lizhou} were about FDFs
of  linear time-invariant and time-varying  deterministic systems,
respectively. Generally speaking, for linear
time-invariant/time-varying  deterministic systems, the
corresopnding FDF design can be turned into solving some LMIs
\cite{liuwangyang,Zhong2003,newtsp1}/differential Riccati equations (DREs)
\cite{liliu,lizhou}. In \cite{Suxiaojie2016}, FDF  of nonlinear
switched stochastic systems was discussed based on T-S fuzzy model
approach.    While  the fault isolation problem  for discrete-time
fuzzy interconnected systems with unknown interconnections  was
considered in \cite{zhzhang}. The reference \cite{liusuntong}
investigated the  adaptive fuzzy output feedback fault-tolerant
optimal control problem for a class of single-input and
single-output nonlinear systems in strict feedback form.  A robust
fault detection $H_-/H_\infty$ observer  was constructed for a
Takagi-Sugeno fuzzy model with sensor faults and unknown bounded
disturbances via an LMI formulation in \cite{Chadli2013}. The
reference \cite{Wangzhenhua2017} studied  $H_-/H_\infty$ fault
detection observer design in the finite frequency domain for a class
of linear parameter-varying  descriptor systems. Recently, we
analyze the $H_-$ index for a class of linear discrete-time
stochastic Markov jump systems with multiplicative noise
\cite{Liyan1}. It can be found that, up to now, there are few works
on mixed $H_-/H_\infty$ FDF for nonlinear It\^o-type stochastic
systems.

Motivated by the aforementioned  reason, this paper studies
$H_-/H_{\infty}$ FDF of affine nonlinear It\^o systems, where a
Luenberger type observer   is considered in designing
$H_-/H_{\infty}$ FDF, which  can not only suppress the effect of
external interference on the residual signal below a level
$\gamma>0$, but also guarantee the residual signal to be sensitive
to the fault signal. The contributions of this paper can be
summarized as follows:

$\bullet$  Applying It\^o  formula together with square completion
technique, some sufficient conditions are obtained for the existence
of  $H_-/H_{\infty}$ FDF of affine stochastic It\^o systems  in
terms of coupled HJIs. As corollaries,  the existence conditions of
$H_-/H_{\infty}$ FDF of linear  stochastic It\^o systems are also
presented in terms of algebraic Riccati inequalities (ARIs), which
can be transformed into solving LMIs by Matlab LMI Toolbox. How to
solve the HJIs is a challenging problem, and a potential powerful
technique to solve the coupled  HJIs can refer to \cite{Abu-Khalaf}
by using the neural network approach.

$\bullet$ We also study $H_-/H_{\infty}$ FDF for a class of
quesi-linear stochastic It\^o systems. It is  shown that for such a
class of nonlinear stochastic systems, the corresponding
$H_-/H_{\infty}$ FDF can be designed via solving some LMIs instead
of HJIs, which is very convenient in practice.

 The organization  of this paper is as follows:   In Section 2, we make some preliminaries to introduce  useful definitions and lemmas. In Section 3,
some sufficient conditions for the existence of $H_-/H_\infty$
fault detection filter are presented  based on HJIs.   Applying
LMI-based technique, in Section 4,  for a class of quasi-linear
stochastic systems, the corresponding $H_-/H_\infty$ FDF design  is converted into solving LMIs.   Section 5 presents an
example  to illustrate the effectiveness of our given results.
Section 6 concludes this paper with some remarks and perspectives.

For convenience,  this paper adopts the following standard
notations:

${\mathcal R}^+:=[0,+\infty)$;  $M'$: the transpose of the matrix
$M$ or vector $M$; $M>0$ ($M<0$): the matrix $M$ is a positive
definite (negative definite) real symmetric matrix; $I_n$: $n\times
n$ identity matrix; ${\mathcal R}^n$: the $n$-dimensional real
Euclidean vector space with the norm
$\|x\|=\sqrt{\sum_{k=0}^{n}x_k^2}$; ${\mathcal R}^{n\times m}$: the
$n\times m$  real  matrix space; $\mathcal {L}^2_{\mathcal
{F}}(\mathcal {R}^+,\mathcal {R}^{n_v})$: the space of
nonanticipative stochastic process $v(t)\in\mathcal {R}^{n_v}$ with
respect to an increasing $\sigma$-algebra $\{\mathcal
{F}_t\}_{t\geq0}$ satisfying $\|v(t)\|_{\mathcal
{L}^2_{\infty}}:=E\int^{\infty}_0\|v(t)\|^2dt<\infty$; A function
$f(x)$ is called a positive function, if $f(x)>0$ for any $x\neq0$,
and $f(0)=0$; $C^2(U;X)$: the class of $X$-valued functions $V(x)$,
which are   twice continuously differential with respect to $x\in
U$, except possibly at the point $x=0$;   $C^{2,1}(U\times {\mathcal
R}^+;X)$: the class of $X$-valued functions $V(x,t)$, which are
twice continuously differential with respect to $x\in U$,  and once
continuously differential with respect to $t\in {\mathcal R}^+$,
except possibly at the point $x=0$.

\section{Preliminaries}

Consider the following affine nonlinear It\^o stochastic system (the time variable $t$ is suppressed):

\begin{flalign}\label{2017fd-system1}
\begin{cases}
dx=\left[f_1(x)+g_1(x)v+h(x)f\right]dt\\
\ \ \ \ \ \ \ \ \ \ \
+\left[f_2(x)+g_2(x)v\right]dw,\\
y=l(x)+m(x)v+n(x)f,\\
x(0)=x_0\in\mathcal {R}^{n},
\end{cases}
\end{flalign}
where  $x(t)\in \mathcal{R}^{n}$ is the $n$-dimensional state
vector, $y\in \mathcal{R}^{n_y}$ is the $n_y$-dimensional
measurement  output, $v\in \mathcal{R}^{n_v}$ stands for the
exogenous disturbance signal with $v\in\mathcal {L}^2_{\mathcal
{F}}(\mathcal {R}^+,\mathcal {R}^{n_v})$, and $w(t)$ is a 1-D standard
Wiener process defined on the complete filtered probability space
$(\Omega,\mathcal {F},\{\mathcal {F}_t\}_{t\in {\mathcal
R}^+},\mathcal {P})$ with the $\sigma$-field  $\mathcal {F}_t$
generated by $w(\cdot)$ up to time $t$. $f(t)\in\mathcal {R}^{n_f}$
denotes the fault information to be detected. We
assume that all functions $f_1(x)$, $g_1(x)$, $h(x)$, $f(t)$,
$f_2(x)$, $g_2(x)$, $l(x)$. $m(x)$ and $n(x)$ are continuous, which
satisfy certain conditions  as linear growth condition and
Lipschitz condition such that the state equation in
(\ref{2017fd-system1}) has a unique strong solution.

In this paper, we adopt the following   Luenberger-type observer as
FDF  of  system (\ref{2017fd-system1}):

\begin{flalign}\label{2017fd-system2}
\begin{cases}
d\hat{x}(t)=\hat{f}(\hat{x}(t))dt+\hat{h}(\hat{x}(t))\left(y(t)-l(\hat{x}(t))\right)dt,\\
r(t)=\hat{s}(\hat{x}(t))(y(t)-l(\hat{x}(t))),\\
\hat{x}(0)=0,
\end{cases}
\end{flalign}
where $\hat{x}(t)$ is  the estimated value of $x(t)$,    $r(t)$ is
viewed as the residual signal,  $\hat{f}$, $\hat{h}$, and $\hat{s}$
are the filter functions to be designed.

Set
$\eta(t)=\left[\begin{array}{cc}x(t)'&\hat{x}(t)'\end{array}\right]'$,
then we get the following augmented system:
\begin{flalign}\label{2017fd-system3}
\begin{cases}
d\eta=\left(\tilde{f}_1(\eta)+\tilde{g}_1(\eta)v+\tilde{h}(\eta)f\right)dt\\
\ \ \ \ \ \ \ \ \ \ \ \
+\left(\tilde{f}_2(\eta)+\tilde{g}_2(\eta)v\right)dw,\\
r(t)=\tilde{s}(\eta,v,f),\\
\eta(0)=\left[\begin{array}{ccc}x_0\\0
\end{array}\right]\in \mathcal {R}^{n_{\eta}},
\end{cases}
\end{flalign}
where
\begin{eqnarray*}
&&\tilde{f}_1(\eta)= \left[\begin{array}{cccc} f_1(x)\\
\hat{f}(\hat{x})+\hat{h}(\hat{x})\left(l(x)-l(\hat{x})\right)
\end{array}\right], \tilde{g}_1(\eta)= \left[\begin{array}{cccc} g_1(x)\\
\hat{h}(\hat{x})m(x(t))
\end{array}\right],\\
&&\tilde{h}(\eta)= \left[\begin{array}{cccc} h(x)\\
\hat{h}(\hat{x})n(x)
\end{array}\right], \tilde{f}_2(\eta)= \left[\begin{array}{cccc} f_2(x)\\0
\end{array}\right],\\
&&\tilde{g}_2(\eta)= \left[\begin{array}{cccc} g_2(x)\\0
\end{array}\right], \tilde{s}(\eta,v,f)=
\hat{s}(\hat{x})\left[l(x)+m(x)v+n(x)f-l(\hat{x})\right].
\end{eqnarray*}
\subsection{Definitions and lemmas}
For our needs, we introduce some definitions as follows.

\begin{definition} An It\^o-type stochastic differential system
\begin{flalign}\label{2017fd-system1ssdd}
\begin{cases}
dx(t)=f(x(t))\,dt+g(x(t))dw(t),\\
x(0)=x_0\in\mathcal {R}^{n}
\end{cases}
\end{flalign}
is said to be exponentially stable in mean square sense, if there
exist constants $\beta\geq1$  and  $\alpha>0$, such that the
solution $x(t)$ of system (\ref{2017fd-system1ssdd}) satisfies
$$
E\|x(t)\|^2\leq\beta e^{-\alpha t}\|x_0\|^2, \ \  t\ge 0.
$$
\end{definition}
The residual signal $r(t)$ needs to be measured and calculated in
real time to judge whether there is a noticeable fault occuring.
Therefore, it is worth emphasizing that the comprehensive ability to
reflect external interferences and internal fault information are
important for $r(t)$ in the designed filter. In order to reasonably
analyze and study  the capacity of $r(t)$, we introduce $H_{-}$
index and $H_{\infty}$ index to discribe system
(\ref{2017fd-system3}).

\begin{definition}
For stochastic system (\ref{2017fd-system3}), its sensitive
operator $\mathcal {L}_{f,r}$ and  $H_{-}$ index are respectively defined as
\begin{eqnarray*}
\mathcal {L}_{f,r}:f(t) \in \mathcal {L}^2_{\mathcal {F}}(\mathcal
{R}^+,\mathcal {R}^{n_f})\mapsto r(t)\in \mathcal {L}^2_{\mathcal
{F}}(\mathcal {R}^+,\mathcal {R}^{n_r})
\end{eqnarray*}
and
\begin{eqnarray*}
\|\mathcal {L}_{f,r}\|_{-}=\inf_{\begin{array}{cc} v(t)\equiv0,
\eta(0)=0,\\ f(t)\not\equiv0,
 f(t)\in \mathcal {L}^2_{\mathcal
{F}}(\mathcal {R}^+,\mathcal {R}^{n_f}) \end{array}}
\frac{\|r(t)\|_{\mathcal {L}^2_{\infty}}}{\|f(t)\|_{\mathcal
{L}^2_{\infty}}}.
\end{eqnarray*}
Meanwhile, we  define the perturbation operator $\mathcal {L}_{v,r}$
and $H_{\infty}$ index as
\begin{eqnarray*}
\mathcal {L}_{v,r}: v(t)\in \mathcal {L}^2_{\mathcal {F}}(\mathcal
{R}^+,\mathcal {R}^{n_v})\mapsto r(t)\in \mathcal {L}^2_{\mathcal
{F}}(\mathcal {R}^+,\mathcal {R}^{n_r})
\end{eqnarray*}
and
\begin{eqnarray*}
\|\mathcal {L}_{v,r}\|_{\infty}=\sup_{\begin{array}{cc}f(t)\equiv0,
\eta(0)=0, \\v(t)\not\equiv0, v(t)\in \mathcal {L}^2_{\mathcal
{F}}(\mathcal {R}^+,\mathcal {R}^{n_v})\end{array}}
\frac{\|r(t)\|_{\mathcal {L}^2_{\infty}}}{\|v(t)\|_{\mathcal
{L}^2_{\infty}}},
\end{eqnarray*}
respectively.
\end{definition}
The purpose   in this paper is to design a mixed $H_{-}/H_{\infty}$
FDF for system (\ref{2017fd-system1}). Below, we
give  the definition of mixed $H_{-}/H_{\infty}$ FDF.

\begin{definition}\label{2017fd-def1}
A FDF (\ref{2017fd-system2}) is called the mixed
$H_{-}/H_{\infty}$ FDF, if for any given scalars
$\gamma>0$ and $\delta>0$, the following requirements are satisfied
simultaneously.
\begin{itemize}
  \item The augmented system (\ref{2017fd-system3}) is internally stable,  that is, when $f(t)\equiv0$ and $v(t)\equiv0$
  in system (\ref{2017fd-system3}), the following system
\begin{flalign}\label{2017fd-systemsss}
\begin{cases}
d\eta(t)=\tilde{f}_1(\eta(t))dt+\tilde{f}_2(\eta(t))\,dw(t),\\
\eta(0)=\left[\begin{array}{ccc}x_0\\0
\end{array}\right]\in \mathcal {R}^{n_{\eta}},
\end{cases}
\end{flalign}
is exponentially stable in mean square  sense.
  \item The augmented system (\ref{2017fd-system3}) is externally stable, that is, when $f(t)\equiv0$,  $\eta(0)=0$ in system
  (\ref{2017fd-system3}) , for any  nonzero  $v(t) \in \mathcal {L}^2_{\mathcal {F}}(\mathcal {R}^+,\mathcal {R}^{n_v})$, the $L_2$-gain from $v(t)$ to $r(t)$
    of the following system
       \begin{flalign}\label{2017fd-system5}
\begin{cases}
d\eta(t)=\left(\tilde{f}_1(\eta(t))+\tilde{g}_1(\eta(t))v(t)\right)dt\\
\ \ \ \ \ \ \ \ \ \ \ \
+\left(\tilde{f}_2(\eta(t))+\tilde{g}_2(\eta(t))v(t)\right)dw(t),\\
r(t)=\hat{s}(\hat{x}(t))\left(l(x(t))+m(x(t))v(t)-l(\hat{x}(t))\right)
\end{cases}
\end{flalign}
       is less than or equal to  $\gamma>0$, i.e., the $H_{\infty}$ index $\|\mathcal {L}_{v,r}\|_{\infty}\leq\gamma$.
  \item The residual signal is enough sensitive to the fault, that is, when $v(t)\equiv0$,  $\eta(0)=0$ in system (\ref{2017fd-system3}),  for any nonzero $f(t)\in \mathcal {L}^2_{\mathcal {F}}(\mathcal {R}^+,\mathcal {R}^{n_f})$, the $L_2$-gain from $f(t)$ to $r(t)$ of the following system
       \begin{flalign}\label{2017fd-system6}
\begin{cases}
d\eta(t)=\left(\tilde{f}_1(\eta(t))+\tilde{h}(\eta(t))f(t)\right)dt+\tilde{f}_2(\eta(t))dw(t),\\
r(t)=\hat{s}(\hat{x}(t))\left(l(x(t))+n(x(t))f(t)-l(\hat{x}(t))\right)
\end{cases}
\end{flalign}
       is large  than or equal to  $\delta>0$, i.e., the $H_{-}$ index $\|\mathcal {L}_{f,r}\|_-\geq\delta$.
\end{itemize}
\end{definition}
In Definition \ref{2017fd-def1},  $\gamma>0$ is called as the disturbance attenuation
level, and $\delta>0$ as the fault sensitivity level.

The following lemmas will play important roles in this study.
\begin{lemma}\cite{book}\label{2017fd-lem1}
For $x, b \in \mathcal {R}^n$, if $A$ is a real symmetric matrix
with appropriate dimension, $A^{-1}$ exists. Then we have
$$
x'Ax+x'b+b'x=(x+A^{-1}b)'A(x+A^{-1}b)-b'A^{-1}b.
$$
\end{lemma}

\begin{lemma}\cite{maobook, book}\label{2017fd-lem22}
Suppose there exists a  function $V(\eta,t)\in C^{2,1}(\mathcal
{R}^{n_{\eta}}\times\mathcal {R}^{+};\mathcal {R})$.  An
infinitesimal generator $\mathcal {L} V (\eta,t): \mathcal
{R}^{n_{\eta}}\times \mathcal {R}^{+} \mapsto \mathcal {R}$
associated with (\ref{2017fd-system3}) is  given by
\begin{eqnarray}\label{2017fd-lv}
&&\mathcal {L}V (\eta,t)\nonumber\\
 &=&\frac{\partial V'}{\partial
t}+\frac{\partial V'}{\partial
\eta}\left(\tilde{f}_1(\eta)+\tilde{g}_1(\eta)v+\tilde{h}(\eta)f(t)\right)
\nonumber\\
&&+\frac{1}{2}\left(\tilde{f}_2(\eta)+\tilde{g}_2(\eta)v\right)'
\frac{\partial^2 V}{\partial
\eta^2}\left(\tilde{f}_2(\eta)+\tilde{g}_2(\eta)v\right),
\end{eqnarray}
and
$$
EV(\eta(t), t)=EV(\eta(t_0), t_0)+E\int_{t_0}^t \mathcal
{L}V(\eta(s), s)\,ds, \ 0\le t_0<t<+\infty.
$$
\end{lemma}

\subsection{Residual evaluation}
 After obtaining the gain matrices of the filter, we are in the position to  discuss the residual estimation.
For the purpose of evaluating the residual signal, one tends to
adopt a threshold $J_{th} > 0$, and the $J_{th}$  conforms to the
following decision logic:
\begin{eqnarray}
\begin{cases}
J_r(t)>J_{th}\Rightarrow \text{faults}\Rightarrow \text{alarm},\\
J_r(t)\leq J_{th}\Rightarrow \text{fault-free},
\end{cases}
\end{eqnarray}
where the residual evaluation function $J_r(t)$ is  defined by
\begin{eqnarray}
J_r(t):=\left(\left\{\frac{1}{t}\int^t_0r'(s)r(s)ds\right\}\right)^{\frac{1}{2}},
J_r(0)=0,
\end{eqnarray}
and the threshold \cite{Lih2017}  is determined by
\begin{eqnarray}
J_{th}:=\sup_{f(t)\equiv0, v(t)\in \mathcal {L}^2_{\mathcal
{F}}(\mathcal {R}^+,\mathcal {R}^{n_v})}E J_r(T),
\end{eqnarray}
where $T$ is the evaluation window.

\section{FDF for  Affine Nonlinear Stochastic Systems}
In this section, we will give our main results about the mixed
$H_-/H_{\infty}$ FDF for affine nonlinear stochastic systems.

\begin{theorem}\label{2017fd-th1}
For any given disturbance attenuation level $\gamma>0$ and fault
sensitivity level $\delta>0$, if there exist positive constants
$c_1$,  $c_2$, $c_3$, $\varepsilon_1$, $\varepsilon_2$ and a
positive Lyapunov function $V(\eta)\in C^2(\mathcal
{R}^{n_{\eta}};\mathcal {R}^+)$, such that
\begin{eqnarray}\label{2017fd-sasa}
c_1\|\eta\|^2\leq V(\eta)\leq c_2\|\eta\|^2,
\end{eqnarray}
and $V(\eta)$ solves the following two coupled HJIs
\begin{eqnarray}
&&\begin{cases}\label{2017fd-inq1}
(1+\varepsilon_1)\|\hat{s}(\hat{x}) \left(l(x)-l(\hat{x})\right)\|^2 +\frac{\partial V'}{\partial \eta}\tilde{f}_1(\eta )\\
\ \ \ \
-\left(\frac{1}{2}\tilde{g}_2(\eta )'\frac{\partial^2 V}{\partial \eta^2}\tilde{f}_2(\eta )+\frac{1}{2}\tilde{g}'_1(\eta )\frac{\partial V}{\partial \eta}\right)'\\
\ \ \ \ \cdot
 \left(\frac{1}{2}\tilde{g}_2(\eta )'\frac{\partial^2
V}{\partial \eta^2}\tilde{g}_2(\eta
)+(1+\varepsilon_1^{-1})\|\hat{s}(\hat{x} )m(x
)\|^2I\right.\\
\ \ \ \ \left. -\gamma^2I\right)^{-1}
 \left(\frac{1}{2}\tilde{g}_2(\eta )'\frac{\partial^2
V}{\partial \eta^2}\tilde{f}_2(\eta )+\frac{1}{2}\tilde{g}'_1(\eta
)\frac{\partial V}{\partial \eta}\right)\\
\ \ \ \ +\frac{1}{2}\tilde{f}_2(\eta )'\frac{\partial^2 V}{\partial
\eta^2}\tilde{f}_2(\eta )+c_3\|\eta\|^2\leq0,
\\
\frac{1}{2}\tilde{g}_2(\eta )'\frac{\partial^2 V}{\partial
\eta^2}\tilde{g}_2(\eta )+(1+\varepsilon_1^{-1})\|\hat{s}(\hat{x}
)m(x )\|^2I\\
\ \ \ \ -\gamma^2I<0,
\end{cases}
\end{eqnarray}
and
\begin{eqnarray}
\begin{cases}\label{2017fd-inq2}
(1-\varepsilon_2)\|\hat{s}(\hat{x} )\left(l(\hat{x} )-l(x
)\right)\|^2-\frac{\partial V'}{\partial \eta}\tilde{f}_1(\eta )\\
\ \ \ \ -\frac{1}{2}\tilde{f}_2(\eta )' \frac{\partial^2 V}{\partial
\eta^2}\tilde{f}_2(\eta )
-\frac{1}{2}\frac{\partial V'}{\partial \eta}\tilde{h}(\eta )\\
\ \ \ \
\cdot\left((1-\varepsilon_2^{-1})\|\hat{s}(\hat{x} )n(x )\|^2I-\delta^2I\right)^{-1}\\
\ \ \ \ \cdot \tilde{h}(\eta )'\frac{\partial V}{\partial \eta}\geq0,\\
(1-\varepsilon_2^{-1})\|\hat{s}(\hat{x} )n(x )\|^2I-\delta^2I>0
\end{cases}
\end{eqnarray}
for some filter functions $\hat{f}(\cdot)$, $\hat{h}(\cdot)$,
$\hat{s}(\cdot)$ with suitable dimensions. Then the desired
$H_-/H_\infty$  FDF is obtained by (\ref{2017fd-system2}).

\end{theorem}

\textbf{Proof}: Firstly, we  show that system (\ref{2017fd-system3})
is internally stable,  or equivalently, the
 system (\ref{2017fd-systemsss})  is exponentially stable in mean square sense.
 For system
(\ref{2017fd-systemsss}),    apply Lemma~\ref{2017fd-lem22} and
consider (\ref{2017fd-inq1}), we have
$$
\mathcal {L}V(\eta)|_{f\equiv0,v\equiv0}=\frac{\partial V'}{\partial
\eta}\tilde{f}_1(\eta
)+\frac{1}{2}\tilde{f}_2(\eta)'\frac{\partial^2 V}{\partial
\eta^2}\tilde{f}_2(\eta)\leq-c_3\|\eta\|^2.
$$
So, by Theorem 4.4 of  \cite{maobook} and (\ref{2017fd-sasa}),  the
 system (\ref{2017fd-systemsss})  is exponentially stable in mean square sense.

Next, we will show that the system (\ref{2017fd-system5}) is also
externally stable, i.e., $\|\mathcal
{L}_{v,r}\|_{\infty}\leq\gamma$.  From Lemma~\ref{2017fd-lem22}, for
any $T>0$,  $v(t) \in  \mathcal {L}^2_{\mathcal {F}}(\mathcal
{R}^+,\mathcal {R}^{n_v})$, and the initial value $\eta(0)=0$, we
have
\begin{eqnarray*}
&&E\int^T_{0}\left(\|r(s)\|^2-\gamma^2\|v(s)\|^2\right) ds\nonumber\\
&&=E\int^T_{0}\left(\|r(s)\|^2-\gamma^2\|v(s)\|^2+\mathcal {L}V(\eta(s))|_{f\equiv0}\right) ds\\
&&\ \ \ \ -EV(\eta(T))+V(\eta(0))\nonumber\\
&&=E\int^T_0\left[\|r(s)\|^2-\gamma^2\|v(s)\|^2+\frac{\partial
V'}{\partial \eta}\left(\tilde{f}_1(\eta(s))\right.\right.
\\
&&\ \ \ \  \left. +\tilde{g}_1(\eta(s))v(s)\right)
\left.+\frac{1}{2}\left(\tilde{f}_2(\eta(s))+\tilde{g}_2(\eta(t))v(s)\right)'\right.\\
&&\ \ \ \ \left. \cdot\frac{\partial^2 V}{\partial
\eta^2}\left(\tilde{f}_2(\eta(t))+\tilde{g}_2(\eta(s))v(s)\right)\right]
\,ds-EV(\eta(T)).\nonumber
\end{eqnarray*}
Note that for any $\varepsilon_1>0$, we have
\begin{eqnarray*}
\|r(s)\|^2\leq(1+\varepsilon_1)\|\hat{s}(\hat{x}(s)\left(l(x(s))-l(\hat{x}(s))\right)\|^2\\
\ \ \ \
+(1+\varepsilon_1^{-1})\|\hat{s}(\hat{x}(s))m(x(s))\|^2\|v(s)\|^2.
\end{eqnarray*}
So
\begin{eqnarray}
&&E\int^T_{0}\left(\|r(s)\|^2-\gamma^2\|v(s)\|^2\right) ds\nonumber\\
&\leq& E\int^T_0\left[(1+\varepsilon_1)\|\hat{s}(\hat{x}(s)\left(l(x(s))-l(\hat{x}(s))\right)\|^2)\right.\nonumber\\
&&\left.+((1+\varepsilon_1^{-1})\|\hat{s}(\hat{x}(s))m(x(s))\|^2-\gamma^2)\|v(s)\|^2 \right.\nonumber\\
&&+\frac{\partial V'}{\partial \eta}\left(\tilde{f}_1(\eta(s))+\tilde{g}_1(\eta(s))v(s)\right)+\frac{1}{2}\left(\tilde{g}_2(\eta(s))v(s)\right)'\nonumber\\
&&\cdot\frac{\partial^2 V}{\partial
\eta^2}\left(\tilde{g}_2(\eta(s))v(s)\right)
+\frac{1}{2}\tilde{f}_2(\eta(s))'\frac{\partial^2 V}{\partial \eta^2}\tilde{f}_2(\eta(s))\nonumber\\
&&\left.+\tilde{f}_2(\eta(s))'\frac{\partial^2 V}{\partial
\eta^2}\tilde{g}_2(\eta(s))v(s)\right] \,ds-EV(\eta(T)).
\end{eqnarray}
By Lemma \ref{2017fd-lem1},
\begin{eqnarray}\label{2017fd-dsa}
&& E\int^T_{0}(\|r(s)\|^2-\gamma^2\|v(s)\|^2) ds\nonumber\\
&\leq& E\int^T_0\left\{(1+\varepsilon_1)\|\hat{s}(\hat{x}(s))\left(l(x(s))-l(\hat{x}(s))\right)\|^2\right. \nonumber\\
&&+\frac{\partial V'}{\partial \eta}\tilde{f}_1(\eta(s))
-\left(\frac{1}{2}\tilde{g}_2(\eta(s))'\frac{\partial^2 V}{\partial
\eta^2}\tilde{f}_2(\eta(s))\right.\nonumber\\
&&\left. +\frac{1}{2}\tilde{g}'_1(\eta(s))\frac{\partial V}{\partial
\eta}\right)'
\Lambda_2^{-1}\left(\frac{1}{2}\tilde{g}_2(\eta(s))'\frac{\partial^2 V}{\partial \eta^2}\tilde{f}_2(\eta(s))\right.\nonumber\\
&&\left.\left. +\frac{1}{2}\tilde{g}'_1(\eta(s))\frac{\partial
V}{\partial
\eta}\right)+\left[v+\Lambda_1\right]'\Lambda_2\left[v+\Lambda_1\right]\right.\nonumber\\
&&\left.+\frac{1}{2}\tilde{f}_2(\eta(s))'\frac{\partial^2
V}{\partial \eta^2}\tilde{f}_2(\eta(s))\right\} ds  -EV(\eta(T)),
\end{eqnarray}
where
\begin{eqnarray*}
&&\Lambda_1=\Lambda_2^{-1}\left(\frac{1}{2}\tilde{g}_2(\eta)'\frac{\partial^2 V}{\partial \eta^2}\tilde{f}_2(\eta)+\frac{1}{2}\tilde{g}'_1(\eta)\frac{\partial V}{\partial \eta}\right),\\
&&\Lambda_2=\frac{1}{2}\tilde{g}_2(\eta)'\frac{\partial^2
V}{\partial
\eta^2}\tilde{g}_2(\eta)+(1+\varepsilon_1^{-1})\|\hat{s}(\hat{x})m(x)\|^2I-\gamma^2I.
\end{eqnarray*}
By (\ref{2017fd-inq1}), $\Lambda_2<0$, then inequality
(\ref{2017fd-dsa}) yields that
\begin{eqnarray}
&& E\int^T_{0}\left(\|r(s)\|^2-\gamma^2\|v(s)\|^2\right) ds\nonumber\\
&\leq& E\int^T_0\left\{(1+\varepsilon_1)\|\hat{s}(\hat{x}(s))\left(l(x(s))-l(\hat{x}(s))\right)\|^2\right.\nonumber\\
 &&+\frac{\partial V'}{\partial \eta}\tilde{f}_1(\eta(s))-\left(\frac{1}{2}\tilde{g}_2(\eta(s))'\frac{\partial^2 V}{\partial
\eta^2}\tilde{f}_2(\eta(s))\right.\nonumber\\
&&\left. +\frac{1}{2}\tilde{g}'_1(\eta(s))\frac{\partial V}{\partial
\eta}\right)'
\Lambda_2^{-1}\left(\frac{1}{2}\tilde{g}_2(\eta(s))'\frac{\partial^2
V}{\partial \eta^2}\tilde{f}_2(\eta(s))\right.\nonumber\\
&&\left.
+\frac{1}{2}\tilde{g}'_1(\eta(s))\frac{\partial V}{\partial \eta}\right)\nonumber\\
&&\left.+\frac{1}{2}\tilde{f}_2(\eta(s))'\frac{\partial^2
V}{\partial \eta^2}\tilde{f}_2(\eta(s))\right\} \,ds-EV(\eta(T)),
\end{eqnarray}
Now, the following inequality  holds under the condition
(\ref{2017fd-inq1}):
$$
E\int^T_{0}(\|r(s)\|^2-\gamma^2\|v(s)\|^2)\,ds\leq0.
$$
Let $T\rightarrow \infty$ in the above inequality, we have
\begin{eqnarray}
&&E\int^{\infty}_{0}\|r(t)\|^2dt \leq \gamma^2 E\int^{\infty}_{0}\|v(t)\|^2\,dt, \  \forall v(t)
\in \mathcal {L}^2_{\mathcal {F}}(\mathcal {R}^+;\mathcal
{R}^{n_v}).  \label{eq ughuh}
\end{eqnarray}
Thus, the external stability of the system (\ref{2017fd-system5}) is
shown.

Finally, we show the  $H_{-}$ index to satisfy  $\|\mathcal
{L}_{f,r}\|_{\infty}\geq\delta$.  By It\^o formula, for any $T>0$,
$t\in[0, T]$, $f(t)\not\equiv0$, $v(t)\equiv0$, and the initial value
$\eta(0)=0$, we have
\begin{eqnarray}
&&E\int^T_{0}\left(\|r(s)\|^2-\delta^2\|f(s)\|^2\right)\,ds\nonumber\\
&&=E\int^T_{0}\left(\|r(s)\|^2-\delta^2\|f(s)\|^2-\mathcal
{L}V(\eta(s))|_{v\equiv
0}\right)ds\nonumber\\
&&\ \ \ \ \ \ \ \ \ \ \ +EV(\eta(T))-V(\eta(0))\nonumber\\
&&=E\int^T_0\left[\|r(s)\|^2-\delta^2\|f(s)\|^2-\frac{\partial
V'}{\partial
\eta}\left(\tilde{f}_1(\eta(s))\right.\right.\nonumber\\
&&\ \ \ \ \left.\left. +\tilde{h}(\eta(s))f(s)\right)
-\frac{1}{2}\tilde{f}_2(\eta(s))' \frac{\partial^2 V}{\partial
\eta^2}\tilde{f}_2(\eta(s))\right]ds\nonumber\\
&&\ \ \ \ \ \ \ \ \ \ \ \ +EV(\eta(T))\nonumber\\
&\geq&E\int^T_0\left[(1-\varepsilon_2)\|\hat{s}(\hat{x}(s))\left(l(\hat{x}(s))-l(x(s))\right)\|^2\right.\nonumber\\
&&\ \ \ \
+(1-\varepsilon_2^{-1})\|\hat{s}(\hat{x}(s))n(x(s))\|^2\|f(s)\|^2
-\delta^2\|f(s)\|^2
\nonumber\\
&&-\frac{\partial V'}{\partial
\eta}\left(\tilde{f}_1(\eta(s))+\tilde{h}(\eta(s))f(s)\right)\nonumber\\
&&\ \ \ \ \left.-\frac{1}{2}\tilde{f}_2(\eta(s))' \frac{\partial^2
V}{\partial \eta^2}\tilde{f}_2(\eta(t))\right]\,ds +EV(\eta(T)).
\end{eqnarray}
The above inequality makes use of the following relation ($\forall
\varepsilon_2>0$)
\begin{eqnarray*}
&&\|\hat{s}(\hat{x})\left(l(x)+n(x)f-l(\hat{x})\right)\|^2\\
&&\geq(1-\varepsilon_2)\|\hat{s}(\hat{x})\left(l(\hat{x})-l(x)\right)\|^2+(1-\varepsilon_2^{-1}) \|\hat{s}(\hat{x})n(x)\|^2\|f\|^2.
\end{eqnarray*}
By Lemma \ref{2017fd-lem1},
\begin{eqnarray}
&&E\int^T_{0}\left(\|r(s)\|^2-\delta^2\|f(s)\|^2\right)ds\nonumber\\
&\geq& E\int^T_0
\left\{(1-\varepsilon_2)\|\hat{s}(\hat{x}(s))\left(l(\hat{x}(s))-l(x(s))\right)\|^2\right.\nonumber\\
&&\left.-c_3\|\eta(s)\|^2-\frac{\partial V'}{\partial
\eta}\tilde{f}_1(\eta(s))-\frac{1}{2}\tilde{f}_2(\eta(s))'
\frac{\partial^2 V}{\partial \eta^2}\tilde{f}_2(\eta(s))\right.\nonumber\\
&&+\left[f(s)-\frac{1}{2}\left(I-\delta^2I\right)^{-1}\tilde{h}(\eta(s))'\frac{\partial V}{\partial \eta}\right]'\nonumber\\
&&\cdot\left((1-\varepsilon_2^{-1})\|\hat{s}(\hat{x}(s))n(x(s))\|^2I-\delta^2I\right)\nonumber\\
&&\cdot\left[f(s)-\frac{1}{2}\left(I-\delta^2I\right)^{-1}\tilde{h}(\eta(s))'\frac{\partial
V}{\partial \eta}\right]\left. -\frac{1}{2}\frac{\partial
V'}{\partial
\eta}\tilde{h}(\eta(s))\right.\nonumber\\
&&\cdot\left((1-\varepsilon_2^{-1})\|\hat{s}(\hat{x}(s))n(x(s))\|^2I
-\delta^2I\right)^{-1}\nonumber\\
&&\left.\cdot\tilde{h}(\eta(s))'\frac{\partial V}{\partial
\eta}\right\}\,ds.
\end{eqnarray}
By  the similar technique used in proving (\ref{eq ughuh}), from
inequality (\ref{2017fd-inq2}), we  obtain that for any  $\forall r(t) \in \mathcal {L}^2_{\mathcal {F}}(\mathcal
{R}^+,\mathcal {R}^{n_r})$,
$$
E\int^{\infty}_{0} \|r(s)\|^2\,ds\geq \delta^2E\int^{\infty}_{0}
\|f(s)\|^2\,ds.
$$
The proof is completed. $\square$

\begin{remark}
In fact, when we consider   $H_{\infty}$ index $\|\mathcal
{L}_{v,r}\|_{\infty}$ and $H_{-}$ index $\|\mathcal
{L}_{f,r}\|_{-}$, the same Lyapunov function is not necessary. That
is, we can  choose two different Lyapunov functions $V_1(\eta)$ and
$V_2(\eta)\in C^2(\mathcal {R}^{n_{\eta}}; {\mathcal R}^+)$ in
calculating
$$
\|r(s)\|^2-\gamma^2\|v(s)\|^2+\mathcal {L}V_1(s)|_{f\equiv0}\le 0
$$
and
$$
\|r(s)\|^2-\delta^2\|f(s)\|^2+\mathcal {L}V_2(s)|_{v\equiv0}\ge 0,
$$
which would reduce the conservatism of Theorem \ref{2017fd-th1}.
\end{remark}

\begin{corollary}\label{2017fd-cor2}
For given disturbance attenuation level $\gamma>0$ and fault
sensitivity level $\delta>0$, if for any $\eta\in \mathcal
{R}^{n_{\eta}}$, there exist positive constants $c_1$, $c_2$, $c_3$,
$\varepsilon_1$, $\varepsilon_2$ and two positive Lyapunov functions
$V_1\in C^2(\mathcal {R}^{n_{\eta}};{\mathcal R}^+)$ satisfying
(\ref{2017fd-sasa}) and  $V_2\in C^2(\mathcal {R}^{n_{\eta}};{\mathcal
R}^+)$, which solve  the following coupled HJIs:
\begin{eqnarray}
&&\begin{cases}\label{2017fd-inq5}
(1+\varepsilon_1)\|\hat{s}(\hat{x}) \left(l(x
)-l(\hat{x})\right)\|^2+c_3\|\eta\|^2+\frac{\partial V_1'} {\partial
\eta}\tilde{f}_1(\eta)\\
\ \ -\left(\frac{1}{2}\tilde{g}_2(\eta)'\frac{\partial^2
V_1}{\partial \eta^2}
\tilde{f}_2(\eta)+\frac{1}{2}\tilde{g}'_1(\eta)\frac{\partial V_1}{\partial \eta}\right)'\\
\ \ \cdot\left(\frac{1}{2}\tilde{g}_2(\eta)'\frac{\partial^2
V_1}{\partial
\eta^2}\tilde{g}_2(\eta)+(1+\varepsilon_1^{-1})\|\hat{s}(\hat{x}
)m(x
)\|^2I\right.\\
\ \ \left.-\gamma^2I\right)^{-1}
\left(\frac{1}{2}\tilde{g}_2(\eta)'\frac{\partial^2 V_1}{\partial
\eta^2}\tilde{f}_2(\eta)+\frac{1}{2}\tilde{g}'_1(\eta)\frac{\partial
V_1}{\partial \eta}\right)\\
\ \ +\frac{1}{2}\tilde{f}_2(\eta)'\frac{\partial^2 V_1}{\partial
\eta^2}\tilde{f}_2(\eta)\leq0,
\\
\frac{1}{2}\tilde{g}_2(\eta)'\frac{\partial^2 V_1}{\partial
\eta^2}\tilde{g}_2(\eta)+(1+\varepsilon_1^{-1})\|\hat{s}(\hat{x}
)m(x )\|^2I\\
\ \ -\gamma^2I<0,
\end{cases}
\end{eqnarray}
and
\begin{eqnarray}
&&\begin{cases}\label{2017fd-inq6}
(1-\varepsilon_2)\|\hat{s}(\hat{x}) \left(l(\hat{x} )-l(x
)\right)\|^2-\frac{\partial V_2'}{\partial
\eta}\tilde{f}_1(\eta)\\
\ \ -\frac{1}{2}\tilde{f}_2(\eta)' \frac{\partial^2 V_2}{\partial
\eta^2}\tilde{f}_2(\eta) -\frac{1}{2}\frac{\partial V'_2}{\partial
\eta}\tilde{h}(\eta)\left((1-\varepsilon_2^{-1})\right.\\
\ \ \left.\cdot\|\hat{s}(\hat{x} )n(x )\|^2I
-\delta^2I\right)^{-1}\tilde{h}(\eta)'\frac{\partial V_2}{\partial \eta}\geq0,\\
(1-\varepsilon_2^{-1})\|\hat{s}(\hat{x} )n(x )\|^2I-\delta^2I>0
\end{cases}
\end{eqnarray}
for some filter functions $\hat{f}$, $\hat{h}$, $\hat{s}$ with
suitable dimensions. Then the desired    $H_-/H_\infty$  FDF
is given   by (\ref{2017fd-system2}).
\end{corollary}

If we consider a special case that $m(x)\equiv0$ and $g_2(x)\equiv0$
in system (\ref{2017fd-system1}),  then, we can get the
following result.

\begin{corollary}
For given disturbance attenuation level $\gamma>0$ and fault
sensitivity level $\delta>0$, if for any $\eta\in \mathcal
{R}^{n_{\eta}}$, there exist constants $c_1$, $c_2$, $c_3$,
$\varepsilon_1$, $\varepsilon_2$ and two positive Lyapunov functions
$V_1\in C^2(\mathcal {R}^{n_{\eta}}; {\mathcal R}^+)$  satisfying
(\ref{2017fd-sasa})  and  $V_2\in C^2(\mathcal {R}^{n_{\eta}}; {\mathcal
R}^+)$  solving the following two coupled  HJIs:
\begin{eqnarray}
&&\label{2017fd-inq5} (1+\varepsilon_1)\|\hat{s}(\hat{x} \left(l(x
)-l(\hat{x} )\right)\|^2+c_3\|\eta\|^2+\frac{\partial V_1'}{\partial
\eta}\tilde{f}_1(\eta) \nonumber\\
&&\ \ \ \ +\frac 1 4 \gamma^{-2}
\left(\tilde{g}'_1(\eta)\frac{\partial V_1}{\partial \eta}\right)'
\left(\tilde{g}'_1(\eta)\frac{\partial
V_1}{\partial \eta}\right)\nonumber\\
&&\ \ \ \ +\frac{1}{2}\tilde{f}_2(\eta)'\frac{\partial^2
V_1}{\partial \eta^2}\tilde{f}_2(\eta)\leq0
\end{eqnarray}
and
\begin{eqnarray}
&&\begin{cases}\label{2017fd-inq6}
(1-\varepsilon_2)\|\hat{s}(\hat{x} \left(l(\hat{x} )-l(x
)\right)\|^2-\frac{\partial V_2'}{\partial
\eta}\tilde{f}_1(\eta)\\
\ \ -\frac{1}{2}\tilde{f}_2(\eta)' \frac{\partial^2 V_2}{\partial
\eta^2}\tilde{f}_2(\eta)
-\frac{1}{2}\frac{\partial V'_2}{\partial \eta}\tilde{h}(\eta)\left((1-\varepsilon_2^{-1})\right.\\
\ \ \left.\cdot\|\hat{s}(\hat{x} )n(x )\|^2I-\delta^2I\right)^{-1}\tilde{h}(\eta)'\frac{\partial V_2}{\partial \eta}\geq0,\\
(1-\varepsilon_2^{-1})\|\hat{s}(\hat{x} )m(x )\|^2I-\delta^2I>0
\end{cases}
\end{eqnarray}
for some filter functions $\hat{f}$, $\hat{h}$, $\hat{s}$ with
suitable dimensions. Then  (\ref{2017fd-system2}) is a
desirable $H_-/H_{\infty}$  FDF  for system
(\ref{2017fd-system1}) when  $m(x)\equiv0$ and $g_2(x)\equiv0$.
\end{corollary}

It is well known that for some practical models, not only the state,
but also the external disturbance \cite{Zhang2006} and fault signal
maybe corrupted by noise. For example, the nonlinear stochastic
$H_{\infty}$ control of It\^o type differential systems with all the
state, control input and external disturbance-dependent noise (
$(x,u,v)$-dependent noise for short) was studied in
\cite{Zhangremark2014}. Therefore, a mixed $H_{-}/H_{\infty}$ FDF  for nonlinear stochastic systems with
$(x,v,f)$-dependent noise deserves further study, which motivates us to  consider
the following system
\begin{flalign}\label{2017fd-system7}
\begin{cases}
dx(t)=\left(f_1(x(t))+g_1(x(t))v(t)+h_1(x(t))f(t)\right)dt\\
\ \ \ \ \ \ \ \ \ \ +\left(f_2(x(t))+g_2(x(t))v(t)+h_2(x(t))f(t)\right)dw(t),\\
y(t)=l(x(t))+m(x(t))v(t)+n(x(t))f(t),\\
x(0)=x_0\in\mathcal {R}^{n}.
\end{cases}
\end{flalign}
As so, we can get the following augmented system:
\begin{flalign}\label{2017fd-system8}
\begin{cases}
d\eta(t)=\left(\tilde{f}_1(\eta(t))+\tilde{g}_1(\eta(t))v(t)+\tilde{h}_1(\eta(t))f(t)\right)dt\\
\ \ \ \ \ \ \ \ \
+\left(\tilde{f}_2(\eta(t))+\tilde{g}_2(\eta(t))v(t)+\tilde{h}_2(\eta(t))f(t)\right)dw(t),\\
r(t)=\tilde{s}(\eta(t)),\\
\eta(0)=\left[\begin{array}{ccc}x_0\\0
\end{array}\right]\in \mathcal {R}^{n_{\eta}},
\end{cases}
\end{flalign}
where
\begin{eqnarray*}
\tilde{h}_1(\eta(t))= \left[\begin{array}{cccc} h_1(x)\\
\hat{h}(\hat{x}(t))n(x(t))
\end{array}\right],\ \
\tilde{h}_2(\eta(t))= \left[\begin{array}{cccc} h_2(x)\\ 0
\end{array}\right].\
\end{eqnarray*}
\begin{theorem}\label{2017fd-th2}
For any  given disturbance attenuation level $\gamma>0$ and fault
sensitivity level $\delta>0$, if for all $\eta\in \mathcal
{R}^{n_{\eta}}$,  the following two coupled  HJIs
\begin{eqnarray}
\begin{cases}\label{2017fd-inq7}
(1+\varepsilon_1)\|\hat{s}(\hat{x}) \left(l(x )-l(\hat{x}
)\right)\|^2+c_3\|\eta\|^2+\frac{\partial V_1'}{\partial
\eta}\tilde{f}_1(\eta)\\
\ \
-\left(\frac{1}{2}\tilde{g}_2(\eta)'\frac{\partial^2 V_1}{\partial \eta^2}\tilde{f}_2(\eta)+\frac{1}{2}\tilde{g}'_1(\eta)\frac{\partial V_1}{\partial \eta}\right)'\\
\ \ \cdot\left(\frac{1}{2}\tilde{g}_2(\eta)'\frac{\partial^2
V_1}{\partial
\eta^2}\tilde{g}_2(\eta)+(1+\varepsilon_1^{-1})\|\hat{s}(\hat{x}
)m(x
)\|^2I\right.\\
\ \ \left. -\gamma^2I\right)^{-1}
\cdot\left(\frac{1}{2}\tilde{g}_2(\eta)'\frac{\partial^2
V_1}{\partial
\eta^2}\tilde{f}_2(\eta)+\frac{1}{2}\tilde{g}'_1(\eta)\frac{\partial
V_1}{\partial \eta}\right)\\
\ \ +\frac{1}{2}\tilde{f}_2(\eta)'\frac{\partial^2 V_1}{\partial
\eta^2}\tilde{f}_2(\eta)\leq0,
\\
\frac{1}{2}\tilde{g}_2(\eta)'\frac{\partial^2 V_1}{\partial
\eta^2}\tilde{g}_2(\eta)+(1+\varepsilon_1^{-1})\|\hat{s}(\hat{x}
)m(x )\|^2I\\
\ \ -\gamma^2I<0
\end{cases}
\end{eqnarray}
and
\begin{eqnarray}
&&\begin{cases}\label{2017fd-inq8}
(1-\varepsilon_2)\|\hat{s}(\hat{x}) \left(l(\hat{x} )-l(x
)\right)\|^2-\frac{\partial V_2'}{\partial
\eta}\tilde{f}_1(\eta)\\
\ \ -\frac{1}{2}\tilde{f}_2(\eta)'
\frac{\partial^2 V_2}{\partial \eta^2}\tilde{f}_2(\eta)-\left(\frac{1}{2}\tilde{h}'_2(\eta)\frac{\partial^2 V_2}{\partial \eta^2}\tilde{f}_2(\eta)\right.\\
\ \ \left.+\frac{1}{2}\tilde{h}_1'(\eta)\frac{\partial V_2}{\partial
\eta}\right)'
\left(-\frac{1}{2}\tilde{h}_2(\eta)'\frac{\partial^2 V}{\partial \eta^2}\tilde{h}_2(\eta)+(1-\varepsilon_2^{-1})\right.\\
\ \ \cdot\left.\|\hat{s}(\hat{x} )n(x
)\|^2I-\delta^2I\right)^{-1}\left(\frac{1}{2}\tilde{h}'_2(\eta)\frac{\partial^2
V_2}{\partial \eta^2}\tilde{f}_2(\eta)\right.\\
\ \ \left.
+\frac{1}{2}\tilde{h}'_1(\eta)\frac{\partial V_2}{\partial \eta}\right)\geq0,\\
-\frac{1}{2}\tilde{h}_2(\eta)'\frac{\partial^2 V_2}{\partial
\eta^2}\tilde{h}_2(\eta)+(1-\varepsilon_2^{-1})\|\hat{s}(\hat{x}
)n(x )\|^2I\\
\ \ -\delta^2I>0,
\end{cases}
\end{eqnarray}
admit a set of  solutions  $(V_1>0$, $V_2>0$, $\hat{f}$, $\hat{h}$, $\hat{s}$,
$c_1>0,c_2>0, c_3>0, \varepsilon_1>0, \varepsilon_2>0)$, where
 $V_1\in C^2(\mathcal {R}^{n_{\eta}};{\mathcal R}^+)$  satisfying (\ref{2017fd-sasa}), $V_2\in C^2(\mathcal {R}^{n_{\eta}};{\mathcal R}^+)$ and $c_1,c_2, c_3, \varepsilon_1, \varepsilon_2$ are positive constants. Then the
FDF   for system (\ref{2017fd-system7}) is
given by  (\ref{2017fd-system2}).
\end{theorem}
\textbf{Proof}: Repeating the same procedure as in Theorem
\ref{2017fd-th1} and Corollary \ref{2017fd-cor2}. The proof is
completed. $\square$

%HJIs (\ref{2017fd-inq7}) and
%(\ref{2017fd-inq8}) are very useful tool to study nonlinear It\^o
%stochastic system.

Below, we consider the following
 linear stochastic system
\begin{flalign}\label{2017fd-system12}
\begin{cases}
dx(t)=\left(A_0x(t)+B_0v(t)+C_0f(t)\right)dt\\
\ \ \ \ \ \ \ \ \ \ +\left(A_1x(t)+B_1v(t)+C_1f(t)\right)dw(t),\\
y(t)=A_2x(t)+B_2v(t)+C_2f(t)
\end{cases}
\end{flalign}
as well as the following linear FDF
\begin{flalign}\label{2017fd-system10}
\begin{cases}
d\hat{x}(t)=\left(\hat{A}\hat{x}(t)+\hat{B}(y(t)-A_2\hat{x}(t))\right)dt,\\
r(t)=\hat{S}(y(t)-A_2\hat{x}(t)).
\end{cases}
\end{flalign}
Set
$\eta(t)=\left[\begin{array}{cc}x(t)'&\hat{x}(t)'
\end{array}\right]'$,
then we get the following augmented system
\begin{flalign}\label{2017fd-system13}
\begin{cases}
d\eta(t)=(\tilde{A}_0\eta(t)+\tilde{B}_0v(t)+\tilde{C}_0f(t))dt\\
\ \ \ \ \ \ \ \ \ \
+(\tilde{A}_1\eta(t)+\tilde{B}_1v(t)+\tilde{C}_1f(t))dw(t),\\
r(t)=\tilde{A}_2\eta(t)+\tilde{B}_2v(t)+\tilde{C}_2f(t),
\end{cases}
\end{flalign}
where
$$
\tilde{A}_0=\left[\begin{array}{cc}A_0&0\\ \hat{B}A_2&\hat{A}-\hat{B}A_2\end{array}\right],\ \
\tilde{B_0}=\left[\begin{array}{cc}B_0\\ \hat{B}B_2\end{array}\right],\ \
\tilde{C_0}=\left[\begin{array}{cc}C_0\\ \hat{B}C_2\end{array}\right],\\
$$
$$
\tilde{A}_1=\left[\begin{array}{cc}A_1&0\\0&0\end{array}\right],\ \
\tilde{B}_1=\left[\begin{array}{cc}B_1\\0\end{array}\right],\ \
\tilde{C}_1=\left[\begin{array}{cc}C_1\\0\end{array}\right],\\
$$
$$
\tilde{A}_2=\left[\begin{array}{cc}\hat{S}A_2&-\hat{S}A_2\end{array}\right],\ \
\tilde{B}_2=\hat{S}B_2,\ \
\tilde{C}_2=\hat{S}C_2.\ \
$$
Using Theorem \ref{2017fd-th2},
it is easy to show the following result.

\begin{theorem}\label{2017fd-th6}
For any  given disturbance attenuation level $\gamma>0$ and fault
sensitivity level $\delta>0$, if there exist two positive
 definite matrices $P_1$ and $P_2$ solving the following ARIs:
{\small
\begin{eqnarray}\label{2017fd-sr1}
\begin{cases}
\mathbb{R}_1:=(1+\varepsilon_1)\varpi+P_1\tilde{A}_0+\tilde{A}_0'P_1-(\tilde{B}_1'P_1\tilde{A}_1+\tilde{B}_0'P_1)'\\
\cdot(\tilde{B}_1'P_1\tilde{B}_1+(1+\varepsilon_1^{-1})\|\hat{S}\tilde{B}_2\|^2I-\gamma^2I)^{-1}\\
\cdot(\tilde{B}_1'P_1\tilde{A}_1+\tilde{B}_0'P_1)+\tilde{A}_1'P_1\tilde{A}_1+cI\leq0,\\
(\tilde{B}_1'P_1\tilde{B}_1+(1+\varepsilon_1^{-1})\|\hat{S}\tilde{B}_2\|^2I-\gamma^2I)<0
\end{cases}
\end{eqnarray}}
and
{\small
\begin{eqnarray}\label{2017fd-sr2}
\begin{cases}
\mathbb{R}_2:=(1-\varepsilon_2)\varpi-P_2\tilde{A}_0-\tilde{A}_0'P_2-(\tilde{C}_1'P_2\tilde{A}_1+\tilde{C}_0'P_2)'\\
\cdot(-\tilde{C}_1'P_2\tilde{C}_1+(1-\varepsilon_2^{-1})\|\hat{S}\tilde{C}_2\|^2I-\delta^2I)^{-1}\\
\cdot(\tilde{C}_1'P_2\tilde{A}_1+\tilde{C}_0'P_2)-\tilde{A}_1'P_2\tilde{A}_1\geq0,\\
(-\tilde{C}_1'P_2\tilde{C}_1+(1-\varepsilon_2^{-1})\|\hat{S}\tilde{C}_2\|^2I-\delta^2I)>0,
\end{cases}
\end{eqnarray}}
where
$$\varpi=\left[\begin{array}{cc}
A_2'\hat{S}'\hat{S}A_2&-A_2'\hat{S}'\hat{S}A_2\\-A_2'\hat{S}'\hat{S}A_2&A_2'\hat{S}'\hat{S}A_2
\end{array}\right]$$
and $c, \varepsilon_1$ and $\varepsilon_2$ are positive constants. Then the
FDF  for system (\ref{2017fd-system12}) is
given by  (\ref{2017fd-system10}).
\end{theorem}
\textbf{Proof}: For system (\ref{2017fd-system13}), we choose two
Lyapunov functions as $V_1(\eta)=\eta'P_1\eta$ and
$V_2(\eta)=\eta'P_2\eta$. By  Theorem \ref{2017fd-th2}, we  obtain
{\small
\begin{eqnarray}\label{fefbh}
\begin{cases}
(1+\varepsilon_1)\|\hat{S}(A_2x-A_2\hat{x})\|^2+c\|\eta\|^2
-(\tilde{B}_1'P_1\tilde{A}_1\eta+\tilde{B}_0'P_1\eta)'\\
(\tilde{B}_1'P_1\tilde{B}_1+(1+\varepsilon_1^{-1})\|\hat{S}
\tilde{B}_2\|^2I-\gamma^2I)^{-1}
\left(\tilde{B}_1'P_1\tilde{A}_1\eta\right.\\
\left.+\tilde{B}_0'P_1\eta\right)
+2\eta'P_1\tilde{A}_0\eta+\eta'\tilde{A}_1'P_1\tilde{A}_1\eta\leq0,\\
\tilde{B}_1'P_1\tilde{B}_1+(1+\varepsilon_1^{-1})\|\hat{S}\tilde{B}_2\|^2I-\gamma^2I<0
\end{cases}
\end{eqnarray}
} and {\small
\begin{eqnarray}\label{ewhvd}
\begin{cases}
(1-\varepsilon_2)\|\hat{S}(A_2x-A_2\hat{x})\|^2-(\tilde{C}_1'P_2\tilde{A}_1\eta+\tilde{C}_0'P_2\eta)'\\
(-\tilde{C}_1'P_2\tilde{C}_1
+(1-\varepsilon_2^{-1})\|\hat{S}\tilde{C}_2\|^2I-\delta^2I)^{-1}
\left(\tilde{C}_1'P_2\tilde{A}_1\eta\right.\\
\left.+\tilde{C}_0'P_2\eta\right)
-2\eta'P_2\tilde{A}_0\eta-\eta'\tilde{A}_1'P_2\tilde{A}_1\eta\geq0,\\
-\tilde{C}_1'P_2\tilde{C}_1+(1-\varepsilon_2^{-1})\|\hat{S}\tilde{C}_2\|^2I-\delta^2I>0.
\end{cases}
\end{eqnarray}}
(\ref{fefbh}) and (\ref{ewhvd})  are equivalent to
\begin{eqnarray}\label{dsds}
\begin{cases}
\eta'\mathbb{R}_1\eta\leq0,\\
\tilde{B}_1'P_1\tilde{B}_1+(1+\varepsilon_1^{-1})\|\hat{S}\tilde{B}_2\|^2I-\gamma^2I<0
\end{cases}
\end{eqnarray}
and
\begin{eqnarray}\label{eqvhvhg}
\begin{cases}
\eta'\mathbb{R}_2\eta\geq0,\\
-\tilde{C}_1'P_2\tilde{C}_1+(1-\varepsilon_2^{-1})\|\hat{S}\tilde{C}_2\|^2I-\delta^2I>0,
\end{cases}
\end{eqnarray}
respectively,  where $\mathbb{R}_1$ and  $\mathbb{R}_2$ are defined
respectively in (\ref{2017fd-sr1}) and (\ref{2017fd-sr2}).
(\ref{dsds}) and (\ref{eqvhvhg}) are equivalent to
 (\ref{2017fd-sr1}) and (\ref{2017fd-sr2}), respectively.
  The proof is thus  completed.
$\square$

Theorem \ref{2017fd-th6} is regarding  $H_-/H_{\infty}$ FDF for
linear stochastic systems based on ARIs. (\ref{2017fd-sr1}) and
(\ref{2017fd-sr2}) can be  converted into LMIs by the technique in
the next section, hence, $H_-/H_\infty$  FDF of the system
(\ref{2017fd-system12}) can be easily designed by Matlab LMI
Toolbox.

\section{LMI-Based Approach for Quasi-Linear Systems}
Generally speaking,  for general nonlinear stochastic systems, it is
not easy to design the $H_-/H_{\infty}$ filter due to the difficulty
in solving  HJIs. However, for a class of special nonlinear
stochastic systems called ``quasi-linear stochastic systems", the
filtering design problem can be converted into solving LMIs  as done
in \cite{zhangfilter}.

We consider   the following quasi-linear stochastic system governed
by It\^o differential equation
\begin{flalign}\label{2017fd-system9}
\begin{cases}
dx(t)=\left(A_0x(t)+F_0(x(t))+B_0v(t)+C_0f(t)\right)dt\\
\ \ \ \ \ \ \ \ \ \ \ +\left(A_1x(t)+F_1(x(t))+B_1v(t)+C_1f(t)\right)dw(t),\\
y(t)=A_2x(t)+B_2v(t)+C_2f(t),
\end{cases}
\end{flalign}
where $F_i(0)=0$, $i=1,2$.   As a matter of fact, in
(\ref{2017fd-system7}), if one takes $g_1(x)=B_0$, $h_1(x)=C_0$,
 $g_2(x)=B_1$, $h_2(x)=C_1$, $l(x)=A_2x$, $m(x)=B_2$, $n(x)=C_2$, and regards $A_0x(t)+F_0(x(t))$ and $A_1x(t)+F_1(x(t))$ as
 the Taylor's series expansions of $f_1(x)$ and $f_2(x)$,  respectively, then the state equation of (\ref{2017fd-system7}) comes
down to the first equation of (\ref{2017fd-system9}).

For the quasi-linear stochastic  system (\ref{2017fd-system9}), we
consider linear FDF (\ref{2017fd-system10}). Set
$\eta(t)=\left[\begin{array}{cc}x(t)'&\hat{x}(t)'\end{array}\right]'$,
then we get the following augmented system:
\begin{flalign}\label{2017fd-system11}
\begin{cases}
d\eta(t)=(\tilde{A}_0\eta(t)+\tilde{F}_0(\eta(t))+\tilde{B}_0v(t)+\tilde{C}_0f(t))dt\\
\ \ \ \ \ \ \ \ \ \
+(\tilde{A}_1\eta(t)+\tilde{B}_1v(t)+\tilde{F}_1(\eta(t))+\tilde{C}_1f(t))dw(t),\\
r(t)=\tilde{A}_2\eta(t)+\tilde{B}_2v(t)+\tilde{C}_2f(t),
\end{cases}
\end{flalign}
where
$$
\tilde{F}_0(\eta(t))=\left[\begin{array}{cc}F_0(x(t))\\0\end{array}\right],\
\
\tilde{F}_1(\eta(t))=\left[\begin{array}{cc}F_1(x(t))\\0\end{array}\right].\\
$$

For any given disturbance attenuation level $\gamma>0$ and fault
sensitivity level $\delta>0$, we want to seek the filtering
parameters $\hat{A}$, $\hat{B}$  and  $\hat{S}$, such that
$\|\mathcal {L}_{v,r}\|_{\infty}<\gamma$ and $\|\mathcal
{L}_{f,r}\|_{-}>\delta$.

\begin{assumption}\label{2017fd-th5}
Suppose there exists a scalar $\alpha>0$, such that
\begin{eqnarray}\label{2017fd-dsadd}
\|\tilde{F}_i(\eta)\|\leq\alpha\|\eta\|,\ \ i=1,2.
\end{eqnarray}
\end{assumption}

\begin{lemma}\label{2017fd-lem2}
Given scalars $\gamma>0$, $\delta>0$, if the following three matrix
inequalities
\begin{eqnarray}\label{2017fd-ccccee}
 0<P\leq\beta I,
\end{eqnarray}
\begin{eqnarray}\label{2017fd-dddd}
&&\left[
\begin{array}{cccccccc}
\tilde{A}_2'\tilde{A}_2+\tilde{A}_0'P+P\tilde{A}_0+P+4\alpha^2\beta I&P\tilde{B}_0&\tilde{A}_1'P&\tilde{A}_1'P&0\\
*&-\gamma^2I+\tilde{B}_2'\tilde{B}_2&0&\tilde{B}_1'P&\tilde{B}_1'P\\
*&*&-P&0&0\\
*&*&*&-P&0\\
*&*&*&*&-P
\end{array}
\right]<0
\end{eqnarray}
and
\begin{eqnarray}\label{2017fd-cccc}
&&\left[
\begin{array}{cccccccc}
\tilde{A}_0'P+P\tilde{A}_0+P+4\alpha^2\beta I-\tilde{A}_2'\tilde{A_2}&P\tilde{C}_0& \tilde{A}_1'P&\tilde{A}_1'P&0\\
*&\delta^2I-\tilde{C}_2'\tilde{C}_2&0&\tilde{C}_1'P&\tilde{C}_1'P\\
*&*&-P&0&0\\
*&*&*&-P&0\\
*&*&*&*&-P
\end{array}
\right]<0
\end{eqnarray}
admit a pair of positive solutions $(P>0, \beta>0)$, then the
augmented system (\ref{2017fd-system11}) is  internally stable.
Moreover, $\|\mathcal {L}_{v,r}\|_{\infty}<\gamma$ and $\|\mathcal
{L}_{f,r}\|_{-}>\delta$.
\end{lemma}
\textbf{Proof}: Firstly, we   prove that system
(\ref{2017fd-system11}) is  internally stable, i.e., the following
system
\begin{equation}\label{eqcgcgvh}
d\eta(t)=(\tilde{A}_0\eta(t)+\tilde{F}_0(\eta(t)))\,dt
+(\tilde{A}_1\eta(t)+\tilde{F}_1(\eta(t)))\,dw(t)
\end{equation}
is exponentially stable in mean square sense. For convenience, in
this section, we denote  $\mathcal {L}_1$ as the infinitesimal
generator  of system (\ref{2017fd-system11}).  To prove the
internal  stability of (\ref{2017fd-system11}),  we take the
Lyapunov function as $V(\eta)=\eta'P\eta$, where $P>0$ is the
solution of (\ref{2017fd-ccccee})-(\ref{2017fd-cccc}). By It\^o
formula, for system (\ref{eqcgcgvh}),
\begin{eqnarray}\label{2017fd-wwww}
&&\mathcal {L}_1V(\eta)|_{v\equiv 0,f\equiv 0}\nonumber\\
&&=\eta'(\tilde{A}_0'P+P\tilde{A}_0+\tilde{A}_1'P\tilde{A}_1)\eta+2\eta'P\tilde{F}_0+\tilde{F}_1'P\tilde{F}_1 +2\eta'\tilde{A}_1'P\tilde{F}_1.
\end{eqnarray}
By (\ref{2017fd-dsadd}), it follows that
\begin{eqnarray*}
&&2\eta'P\tilde{F}_0+\tilde{F}_1P\tilde{F}_1+2\eta'\tilde{A}_1'P\tilde{F}_1\leq\eta'P\eta+\tilde{F}_0'P\tilde{F}_0\\
&&+2\tilde{F}_1P\tilde{F}_1+\eta'\tilde{A}_1'P\tilde{A}_1\eta\leq\eta'(P+\tilde{A}_1'P\tilde{A}_1+3\alpha^2\beta
I)\eta.
\end{eqnarray*}
Substituting the above inequality into (\ref{2017fd-wwww}), it
yields that
\begin{eqnarray}
\mathcal {L}_1V(\eta)|_{v\equiv 0,f\equiv
0}\leq\eta'(\tilde{A}_0'P+P\tilde{A}_0+2\tilde{A}_1'P\tilde{A}_1+P+3\alpha^2\beta
I )\eta.
\end{eqnarray}
Obviously, under the condition (\ref{2017fd-dddd}), we have
$$
-\theta:
=\lambda_{max}(\tilde{A}_0'P+P\tilde{A}_0+2\tilde{A}_1'P\tilde{A}_1+P+3\alpha^2\beta
I)<0.
$$
So
$$
\mathcal {L}_1V(\eta)|_{v\equiv 0,f\equiv 0} \leq-\theta\|\eta\|^2.
$$
According to \cite{maobook}, system (\ref{2017fd-system11}) is
exponentially stable in mean square sense.

Secondly, we  prove $\|\mathcal {L}_{v,r}\|_{\infty}<\gamma$ with
$\eta(0)=0$ and $f(t)\equiv0$ in (\ref{2017fd-system11}). Note that
when $f(t)\equiv0$  and $\eta(0)=0$, (\ref{2017fd-system11}) becomes
\begin{flalign}\label{2017fd-system11aa}
\begin{cases}
d\eta(t)=(\tilde{A}_0\eta(t)+\tilde{F}_0(\eta(t))+\tilde{B}_0v(t))\,dt\\
\ \ \ \ \ \ \  \ \
+(\tilde{A}_1\eta(t)+\tilde{B}_1v(t)+\tilde{F}_1(\eta(t)))\,dw(t),\\
\eta(0)=0,\\
r(t)=\tilde{A}_2\eta(t)+\tilde{B}_2v(t).
\end{cases}
\end{flalign}
For system (\ref{2017fd-system11aa}), using It\^o formula to
$V(\eta)=\eta'P\eta$, we have
\begin{eqnarray*}
&&E\int_{0}^{T}(\|r(t)\|^2-\gamma^2\|v(t)\|^2)\,dt\nonumber\\
&&=E\int_{0}^{T}[\mathcal {L}_1V(\eta(t))|_{f\equiv
0}+\|r(t)\|^2-\gamma^2\|v(t)\|^2 ] dt\nonumber\\
&&\ \ \ \ \ \ \ \ \ \ \ \ -EV(\eta(T))\nonumber\\
&&\leq
E\int_{0}^{T}\left[\eta'(t)\left(\tilde{A}_2'\tilde{A_2}+\tilde{A}_0'P+P\tilde{A}_0+2\tilde{A}_1'P\tilde{A}_1+P\right.\right.\nonumber\\
&&\ \ \left. \left.+3\alpha^2\beta I\right)\eta(t)
+2\eta(t)'P\tilde{B}_0v(t)+v(t)'\tilde{B}_1'P\tilde{B}_1v(t)\right.\nonumber\\
&&\ \
\left.+2\eta(t)'\tilde{A}_1'P\tilde{B}_1v(t)+2\tilde{F}_1'P\tilde{B}_1v(t)
+v(t)'\right.\nonumber\\
&&\ \ \left.\cdot(-\gamma^2I+\tilde{B}_2'\tilde{B}_2)v(t)
\right]dt-EV(\eta(T))\nonumber\\
&&\leq
E\int_{0}^{T}\left[\eta'\left(\tilde{A}_2'\tilde{A_2}+\tilde{A}_0'P+P\tilde{A}_0+2\tilde{A}_1'P\tilde{A}_1
+P\right.\right.\nonumber\\
&&\ \ \left. \left.+4\alpha^2\beta I\right)\eta
+2\eta(t)'P\tilde{B}_0v(t)+2v(t)'\tilde{B}_1'P\tilde{B}_1v(t)\right.\nonumber\\
&&\ \
\left.+2\eta(t)'\tilde{A}_1'P\tilde{B}_1v(t)+v(t)'(-\gamma^2I+\tilde{B}_2'\tilde{B}_2)v(t)
\right]\,dt\nonumber\\
&&\ \ \ \ \ \ \ \ \ \ \ -EV(\eta(T))\nonumber\\
&&=E\int_{0}^{T} \left[\begin{array}{cc}\eta(t)'&v(t)'
\end{array}\right]\prod\left[\begin{array}{cc}\eta(t)\\v(t)
\end{array}\right]\,dt-EV(\eta(T))\nonumber\\
&&\le  E\int_{0}^{T} \left[\begin{array}{cc}\eta(t)'&v(t)'
\end{array}\right]\prod\left[\begin{array}{cc}\eta(t)\\v(t)
\end{array}\right]\,dt,
\end{eqnarray*}
where
\begin{eqnarray*}
\prod:=\left[\begin{array}{cc} \Pi_{11}
&P\tilde{B}_0+\tilde{A}_1'P\tilde{B}_1\\
*&-\gamma^2I+\tilde{B}_2'\tilde{B}_2+2\tilde{B}_1'P\tilde{B}_1
\end{array}\right]
\end{eqnarray*}
with $$ \Pi_{11}=\tilde{A}_2'\tilde{A_2}+
\tilde{A}_0'P+P\tilde{A}_0+2\tilde{A}_1'P\tilde{A}_1
+P+4\alpha^2\beta I.
$$
By Schur's  complement and inequality (\ref{2017fd-dddd}), we get
$\prod<0$. Thus, for any $T>0$,
\begin{eqnarray*}
&&E\int_{0}^{T}\|r(t)\|^2\,dt \\
&&\le \gamma^2 E\int_{0}^{T} \|v(t)\|^2\,dt+\lambda_{max}(\prod)E\int_{0}^{T} \|{\mathcal Z}(t)\|^2\,dt\\
&&\le (\gamma^2+\lambda_{max}(\prod))E\int_{0}^{T}\|v(t)\|^2\,dt,
\end{eqnarray*}
where ${\mathcal Z}(t)=\left[\begin{array}{cc}\eta(t)'&v(t)'
\end{array}\right]'$.
Let $T\to\infty$ in above, then for any nonzero  $v\in \mathcal
{L}^2_{\mathcal {F}}(\mathcal {R}^+,\mathcal {R}^{n_v})$,
$\|\mathcal {L}_{v,r}\|_{\infty}\le
(\gamma^2+\lambda_{min}(\prod))^{1/2}<\gamma$, so  the external
stability of  (\ref{2017fd-system11}) is proved.

Finally, we show $\|\mathcal {L}_{f,r}\|_{-}>\delta$ when $v\equiv
0$ and $\eta(0)=0$ in  (\ref{2017fd-system11}). Consider  system
\begin{flalign}\label{2017fd-system11fgvf}
\begin{cases}
d\eta(t)=(\tilde{A}_0\eta(t)+\tilde{F}_0(\eta(t))+\tilde{C}_0f(t))dt\\
\ \ \ \ \ \ \ \ \ \ \ \
+(\tilde{A}_1\eta(t)+\tilde{F}_1(\eta(t))+\tilde{C}_1f(t))dw,\\
\eta(0)=0,\\
r(t)=\tilde{A}_2\eta(t)+\tilde{C}_2f(t).
\end{cases}
\end{flalign}
According to Lemma~\ref{2017fd-lem22}, for system
(\ref{2017fd-system11fgvf}),   we are in a position to show  that
\begin{eqnarray}
&&E\int_{0}^{T}(\|r(t)\|^2-\delta^2\|f(t)\|^2)\,dt\nonumber\\
&&=E\int_{0}^{T}[\|r(t)\|^2-\delta^2\|f(t)\|^2-\mathcal
{L}_1V(\eta(t))|_{v\equiv 0}]\,dt\nonumber\\
&&\ \ \ \ \ \ \ \
+EV(\eta(t))-EV(\eta(0))\nonumber\\
&&\geq
E\int_{0}^{T}\left[-\eta'(t)\left(-\tilde{A}_2'\tilde{A_2}+\tilde{A}_0'P+P\tilde{A}_0+2\tilde{A}_1'P\tilde{A}_1\right.\right.\nonumber\\
&&\ \ \left.\left.+P +4\alpha^2\beta I\right)\eta(t)
-2\eta(t)'P\tilde{C}_0f(t)\right.\nonumber\\
&&\left.\ \ -2f(t)'\tilde{C}_1'P\tilde{C}_1f(t)
-2\eta(t)'\tilde{A}_1'P\tilde{C}_1f(t)\right.\nonumber\\
&&\ \ \left.+f(t)'(-\delta^2I+\tilde{C}_2'\tilde{C}_2)f(t)
\right]dt+EV(\eta(T))\nonumber\\
&&=E\int_{0}^{T} \left[\begin{array}{cc}\eta(t)'&f(t)'
\end{array}\right]\bigsqcup\left[\begin{array}{cc}\eta(t)\\f(t)
\end{array}\right]dt+EV(\eta(T)),\end{eqnarray}
where, by (\ref{2017fd-cccc}),
\begin{eqnarray*}
\bigsqcup:=\left[\begin{array}{cc} \Phi_{11}&
-P\tilde{C}_0-\tilde{A}_1'P\tilde{C}_1
\\
*&-\delta^2I+\tilde{C}_2'\tilde{C}_2-2\tilde{C}_1'P\tilde{C}_1
\end{array}\right]>0
\end{eqnarray*}
with
$$
 \Phi_{11}=\tilde{A}_2'\tilde{A_2}-\tilde{A}_0'P-P\tilde{A}_0-2\tilde{A}_1'P\tilde{A}_1-P-4\alpha^2\beta
I.
 $$
 So
 {\small
\begin{eqnarray}\label{2017fd-das2}
&&E\int_{0}^{T}\|r(t)\|^2\,dt\nonumber\\
&&\ge  \delta^2
E\int_{0}^{T}\|f(t)\|^2\,dt+\lambda_{\min}(\bigsqcup)E\int_{0}^{T}
\|{\mathcal M}(t)\|^2\,dt, \
\end{eqnarray}
} where ${\mathcal M}(t)=\left[\begin{array}{cc}\eta(t)'&f(t)'
\end{array}\right]'$.
For $f(t)\not\equiv 0$, $f\in \mathcal {L}^2_{\mathcal {F}}(\mathcal
{R}^+,\mathcal {R}^{n_f})$,  let $T\to\infty$ in
(\ref{2017fd-das2}), we have
\begin{eqnarray*}\label{2017fd-das2ss}
&&E\int_{0}^{\infty}\|r(t)\|^2\,dt\\
&&\ge  \delta^2 E\int_{0}^{\infty}\|f(t)\|^2\,dt+\lambda_{\min}(\bigsqcup)E\int_{0}^{\infty} \|{\mathcal M}(t)\|^2\,dt\\
&&\ge  \delta^2
E\int_{0}^{\infty}\|f(t)\|^2\,dt+\lambda_{\min}(\bigsqcup)E\int_{0}^{\infty}
\|f(t)\|^2\,dt
\end{eqnarray*}
which yields that $\|\mathcal {L}_{f,r}\|_{-}\ge
(\delta^2+\lambda_{\min}(\bigsqcup))^{1/2}>\delta$ due to
$\lambda_{\min}(\bigsqcup)>0$. The proof is completed.  $\square$

Based on Lemma \ref{2017fd-th5}, an $H_-/H_{\infty}$  FDF  can be
designed in terms of LMIs for system (\ref{2017fd-system9}). For
convenience, we  set   $P$  a  real symmetric diagonal matrix such
as $P=\left[\begin{array}{cccc} P_1&0\\ 0 &P_2
\end{array}\right]>0$, then
\begin{eqnarray*}\label{2017fd-3fsd}
&&\mathcal {A}_0:=P\tilde{A}_0=\left[
\begin{array}{cccccccc}
P_1A_0&0\\
\check{B}A_2&\check{A}-\check{B}A_2
\end{array}
\right], \\
&& \mathcal {B}_0:=P\tilde{B}_0=\left[
\begin{array}{cccccccc}
P_1B_0\\ \check{B}B_2
\end{array}
\right],\\
&& \mathcal {C}_0:=P\tilde{C}_0=\left[
\begin{array}{cccccccc}
P_1C_0\\ \check{B}C_2
\end{array}
\right],\ \ \mathcal {A}_1:=P\tilde{A}_1=\left[
\begin{array}{cccccccc}
P_1A_1&0\\ 0&0
\end{array}
\right],\\
&& \mathcal {B}_1:=P\tilde{B}_1=\left[
\begin{array}{cccccccc}
P_1B_1\\ 0
\end{array}
\right],\ \ \mathcal {C}_1:=P\tilde{C}_1=\left[
\begin{array}{cccccccc}
P_1C_1\\ 0
\end{array}
\right],\ \ \\
&&\mathcal {A}_2:=\tilde{A}_2'\tilde{A}_2=\left[
\begin{array}{cccccccc}
A_2'\check{S}A_2&-A_2'\check{S}A_2\\-A_2'\check{S}A_2&A_2'\check{S}A_2
\end{array}
\right],\\
&&\mathcal {B}_2:=\tilde{B}_2'\tilde{B}_2=B_2'\check{S}B_2,\ \
\mathcal {C}_2:=\tilde{C}_2'\tilde{C}_2=C_2'\check{S}C_2.
\end{eqnarray*}
where $\left[
\begin{array}{cccccccc}
\check{A}&\check{B}
\end{array}
\right]=P_2\left[
\begin{array}{cccccccc}
\hat{A}&\hat{B}
\end{array}
\right]$, $\check{S}=\hat{S}'\hat{S}$.

\begin{theorem}\label{2017fd-th5}
If there exists the  solution $(P_1>0, P_2>0, \beta>0, \check{A},
\check{B}, \check{S})$ solving the following LMIs:
\begin{eqnarray}
0<\left[
\begin{array}{cccccccc}
P_1&0\\*&P_2
\end{array}
\right]\leq \beta I,
\end{eqnarray}
\begin{eqnarray}
\left[
\begin{array}{cccccccc}
\hbar_{11}&\mathcal {B}_0&\mathcal {A}_1'&\mathcal {A}_1'&0\\
*&-\gamma^2I+\mathcal {B}_2&0&\mathcal {B}_1'&\mathcal {B}_1'\\
*&*&-P&0&0\\
*&*&*&-P&0\\
*&*&*&*&-P&\\
\end{array}
\right]<0
\end{eqnarray}
and
\begin{eqnarray}
\left[
\begin{array}{cccccccc}
\ell_{11}&\mathcal {C}_0&\mathcal {A}_1'&\mathcal {A}_1'&0\\
*&\delta^2I-\mathcal {C}_2&0&\mathcal {C}_1'&\mathcal {C}_1'\\
*&*&-P&0&0\\
*&*&*&-P&0\\
*&*&*&*&-P
\end{array}
\right]<0,
\end{eqnarray}
with
$$
\hbar_{11}=\mathcal {A}_2+\mathcal {A}_0'+\mathcal
{A}_0+P+4\alpha^2\beta I
$$
and
$$
\ell_{11}=-\mathcal {A}_2+\mathcal {A}_0'+\mathcal
{A}_0+P+4\alpha^2\beta I,
$$
then  (\ref{2017fd-system10}) is a mixed $H_-/H_{\infty}$ FDF  of
the system (\ref{2017fd-system9}). In this case, the admissible
filter matrices can be given by
\begin{eqnarray*}
\hat{S}=\check{S}^{\frac{1}{2}},\ \ \hat{A}=P_2^{-1}\check{A}, \ \
\hat{B}=P_2^{-1}\check{B}.
\end{eqnarray*}
\end{theorem}

\section{Numerical Example}
In this section, one numerical example is provided to illustrate the
effectiveness of our main results.

\begin{example}
Consider the nonlinear stochastic system (\ref{2017fd-system9}) with
the following parameters:
\begin{eqnarray*}
&&A_0=\left[\begin{array}{ccc} -6.01&-2.94\\-2.94&-6.17
\end{array}\right], \ \
A_1=\left[\begin{array}{ccc} 0.91&-0.44\\-1.31&-0.39
\end{array}\right],\\
&&A_2=\left[\begin{array}{ccc} 0.43&0.15\\-0.09&0.07
\end{array}\right],\ \
B_0=\left[\begin{array}{ccc} 1.37\\-0.41
\end{array}\right],\\
&&B_2=\left[\begin{array}{ccc} 0.35\\-0.6
\end{array}\right],\ \
C_0=\left[\begin{array}{ccc} 1.21\\-0.11
\end{array}\right],\ \
C_2=\left[\begin{array}{ccc} -3.67\\0.51
\end{array}\right],\\
&&F_1(x(t))=B_1=C_1=\left[\begin{array}{ccc} 0\\0
\end{array}\right],\\
&&F_0(x(t))=0.5\left[\begin{array}{ccc} sin(x_1(t))\\sin(x_2(t))
\end{array}\right].
\end{eqnarray*}
In addition,  we choose the $H_\infty$ performance level $\gamma=1$
and $H_-$ performance level  $\delta=0.5$. For the above parameters,
by using Matlab LMI Toolbox, the solutions of LMIs in Theorem
\ref{2017fd-th5} for $ \left\{P_1>0, P_2>0, \beta>0\right.$,
$\left.\check{A}, \check{B}, \check{S}\right\}$ are obtained as
\begin{eqnarray*}
&&P_1=\left[\begin{array}{ccc}
 0.5248&-0.0397\\
   -0.0397&0.3805
\end{array}\right],\ \
P_2=\left[\begin{array}{ccc}
 0.4799&0\\
    0&0.4799
\end{array}\right],\\
&& \check{A}=\left[\begin{array}{ccc}
   -2.0747&-0.8242\\
    0.8209&-2.0835
\end{array}\right],\ \
\check{B}=\left[\begin{array}{ccc}
    0.0032&-0.0529\\
   -0.0054&-0.0441
\end{array}\right],\\
&&\check{S}=\left[\begin{array}{ccc}
    0.4052&0.2562\\
    0.2562&0.2745
\end{array}\right],\ \
\beta=6.
\end{eqnarray*}
Thus, the desired filter matrices  of the $H_-/H_{\infty}$ FDF are
as follows:
\begin{eqnarray*}
&&\hat{A}=\left[\begin{array}{ccc}
     -4.3228&-1.7172\\
    1.7103&-4.3411
\end{array}\right],\ \
\hat{B}=\left[\begin{array}{ccc}
     0.0067&-0.1102\\
   -0.0113&-0.0919
\end{array}\right],\\
&&\hat{S}=\left[\begin{array}{ccc}
     0.0067&-0.1102\\
   -0.0113&-0.0919
\end{array}\right].
\end{eqnarray*}
We use Matlab to simulate the state trajectories $x(t)$ and the
filter trajectories ${\hat x}(t)$  of system (\ref{2017fd-system11})
under $v(t)\equiv0$,  $f(t)\equiv0$ and  $
x_0=\left[\begin{array}{ccc} 1\\-1
\end{array}\right]$; see  Figures \ref{x1x2} and \ref{hatx}. From Figures~\ref{x1x2}-\ref{hatx}, it is easy to see that
system (\ref{2017fd-system11})  is internally stable.

\begin{figure}[!htb]
  \centering
  \includegraphics[height=6cm]{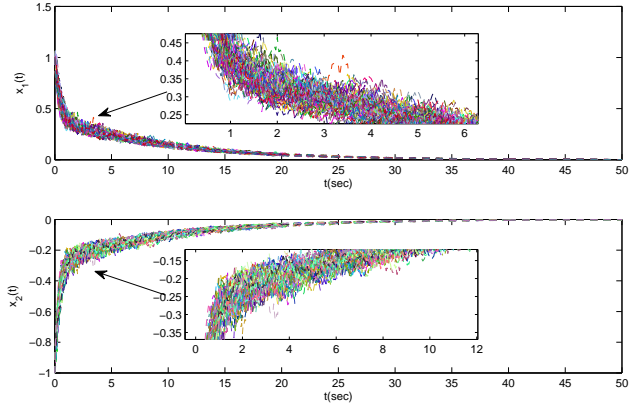}
  \caption{State trajectories $x(t)$ of the
system  (\ref{2017fd-system9}) with $f(t)\equiv0,v(t)\equiv0$.}
  \label{x1x2}
\end{figure}
\begin{figure}[!htb]
  \centering
  \includegraphics[height=6cm]{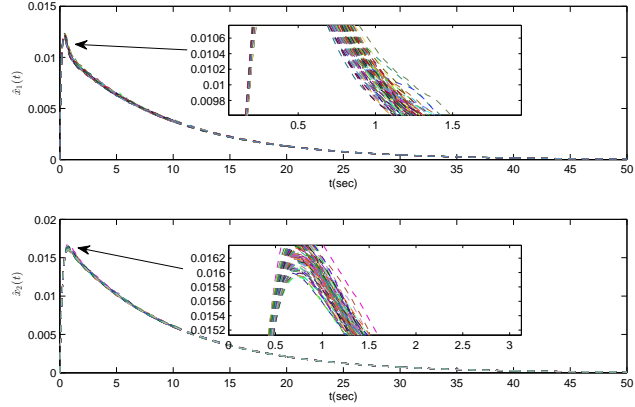}
  \caption{State trajectories $\hat{x}(t)$ of the
system  (\ref{2017fd-system10}).}
  \label{hatx}
\end{figure}
\begin{figure}[!htb]
  \centering
  \includegraphics[height=6cm]{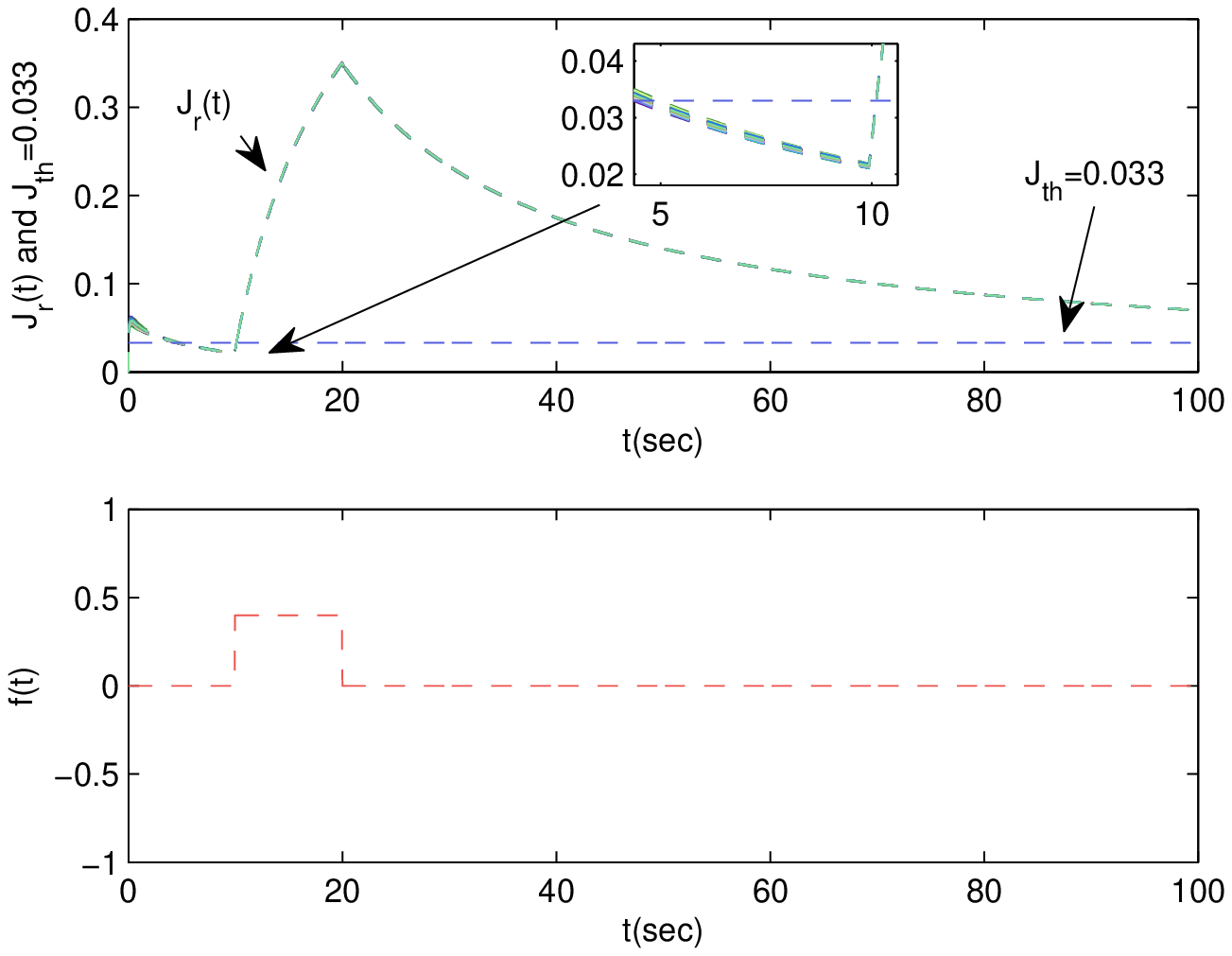}
  \caption{$J_r(t)$, $J_{th}$ and $f(t)=f_{0.4}(t)$.}
  \label{f=0.4}
\end{figure}

\begin{figure}[!htb]
  \centering
  \includegraphics[height=6cm]{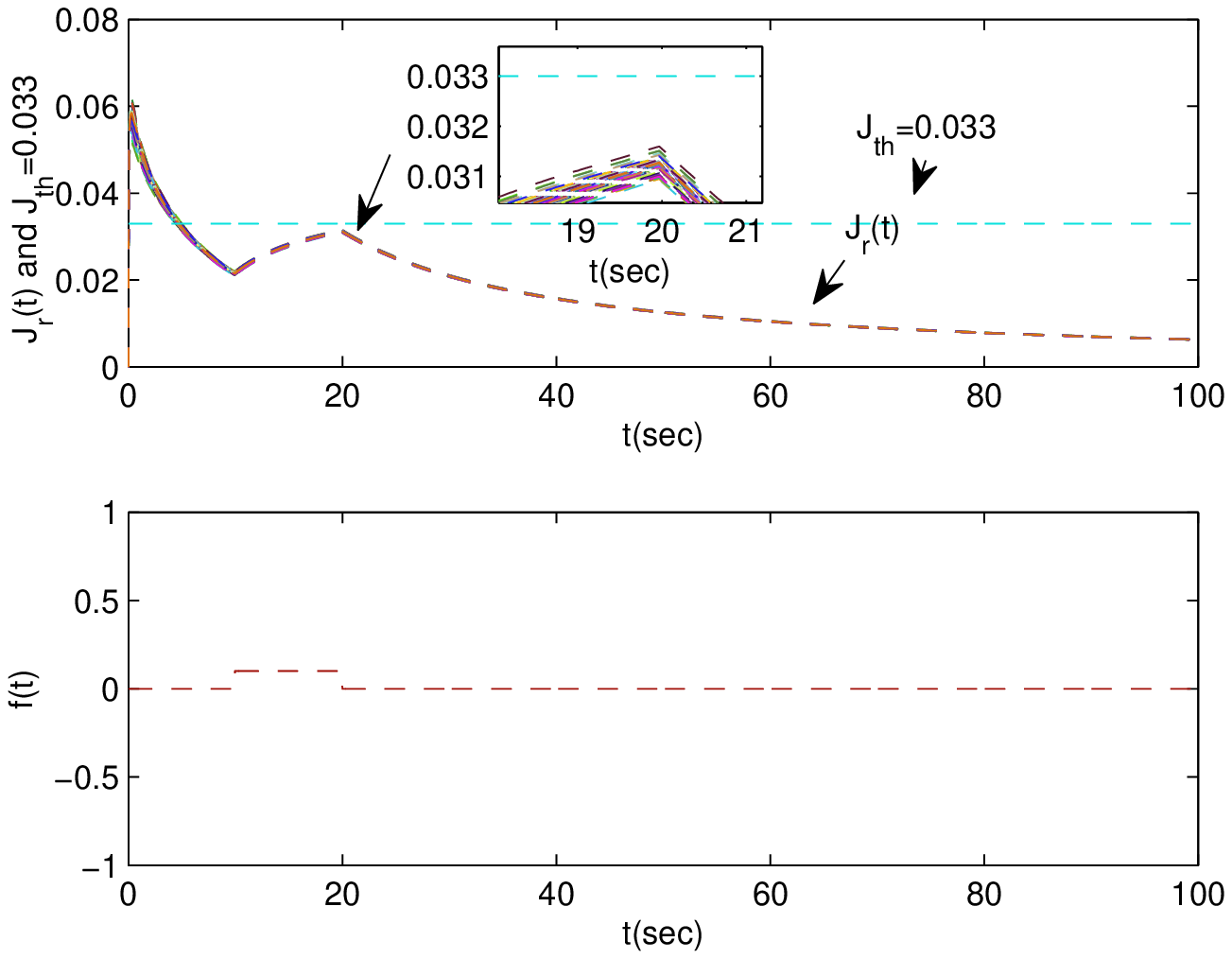}
  \caption{$J_r(t)$, $J_{th}$ and $f(t)=f_{0.1}(t)$.}
  \label{f=0.1}
\end{figure}

\begin{figure}[!htb]
  \centering
  \includegraphics[height=6cm]{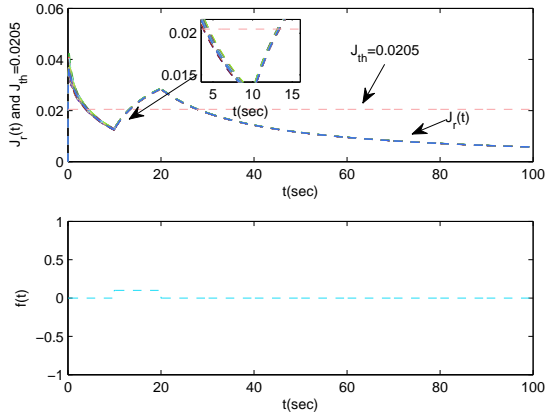}
  \caption{$J_r(t)$, $J_{th}$ and $f(t)=f_{0.1}(t)$  for  $\gamma=0.78$ and  $\delta=1.41$.}
  \label{newf=0.1}
\end{figure}
To show the effectiveness of the designed  $H_-/H_\infty$  filter,
we assume $v(t)=0.9^{t}\in\mathcal {L}^2_{\mathcal {F}}(\mathcal {R}^+,\mathcal {R})$ and $f(t)=f_{0.4}(t)=\begin{cases}0.4,\  \ t\in[10,20]\\
0, \  \ else\end{cases}\in\mathcal {L}^2_{\mathcal {F}}(\mathcal
{R}^+,\mathcal {R})$.   The residual evaluation function is
$$
J_r(t)=\left(E\left\{\frac{1}{t}\int^t_0r'(s)r(s)ds\right\}\right)^{\frac{1}{2}}.
$$
After 100 times Monte Carlo simulations without fault influence,
$J_{th}=0.4$ with evaluation window $T=5$.
The residual evaluation function $J_r(t)$ and fault signal $f(t)$ are depicted in Figure \ref{f=0.4}. However, if the fault
signal is weaker, the FDF may fail to alarm, which  can be seen from Figure \ref{f=0.1}
when we set $f(t)=f_{0.1}(t)=\begin{cases}0.1,\  \ t\in[10,20]\\
0,\ \ \ \ else\end{cases}\in\mathcal {L}^2_{\mathcal {F}}(\mathcal
{R}^+,\mathcal {R})$. In practice, one needs to     select a
suitable combination of $\gamma$ and $\delta$ according to practical
engineering  requirements.  In this example,  if we set
$\gamma=0.78$,  $\delta=1.41$, $f(t)=f_{0.1}(t)$,  the
$H_-/H_{\infty}$ FDF parameters   can be computed  as
\begin{eqnarray*}
&&\hat{A}=\left[\begin{array}{ccc}
     -4.8359&-3.3771\\
    3.3761&-4.509
\end{array}\right],\ \
\hat{B}=\left[\begin{array}{ccc}
     -0.0177&-0.1615\\
   -0.0124&-0.1198
\end{array}\right],\\
&&\hat{S}=\left[\begin{array}{ccc}
     0.6030&0.2787\\
   0.2787&0.4084
\end{array}\right]
\end{eqnarray*}
as well as $J_{th}=0.0205$. From Figure \ref{newf=0.1},  we can see
that  the fault sensitivity is  improved.

\end{example}

\section{Conclusion}
In this paper, the $H_-/H_{\infty}$ FDF design for It\^o-type affine
nonlinear stochastic systems and quasi-linear systems have been
discussed, and sufficient conditions for the existence of the
desired FDF have been  given via HJIs and LMIs, respectively.  The
key to  the FDF design of affine nonlinear stochastic systems is how
to solve the coupled HJIs, this is a very challenging problem, and
some potential effective approaches to overcome this difficulty may
refer to \cite{bschen2009,Deng}
 for
global linearization method, \cite{Abu-Khalaf} for
neural network
method, and \cite{chenchen} for fuzzy interpolation
method. In
addition, from our simulation example, in order to
 select a more suitable combination of $(\gamma,\delta)$, a
co-design algorithm for $H_-$ index and $H_{\infty}$ index is
necessary. Wherever possible, the smaller $\gamma>0$ and the larger
$\delta>0$, the better the performance of FDF.
 However, from the second inequalities of
  HJIs (\ref{2017fd-inq1}) and (\ref{2017fd-inq2}),
   $\gamma>0$ cannot be arbitrarily small and $\delta>0$
   cannot be arbitrarily large. In order to obtain the
   optimal selection  $(\gamma*,\delta*)$, we may turn to Pareto optimization method \cite{yaning1,yaning3} together with convex optimization \cite{Boydbook}. Pareto optimization is a co-operative game, its application to $H_-/H_{\infty}$ FDF design will be our future work.


\begin{thebibliography}{99}


\bibitem{Abu-Khalaf}
 M. Abu-Khalaf, J. Huang, and F. L. Lewis,
 \textit{ Nonlinear $H_2/H_\infty$ Constrained Feedback Control:
   A Practical Design Approach Using Neural Networks}.
   Springer-Verlag: London, 2006.

%%%%$$$$$
\bibitem{berman_06_1}
N. Berman and U. Shaked,  ``$H_\infty$-like control for nonlinear
stochastic control,''
 {\em Syst. Control Lett.},
vol. 55, pp. 247--257, 2006.

%%%%%%%%%$$$$$$$$$
\bibitem{Boydbook} S. Boyd, L. Vandenberghe,
\textit{Convex Optimization}.
 Cambridge, UK: Cambridge University Press, 2004.

%%%%%%%%%%%%$$$$$$
\bibitem{Chadli2013}
M. Chadli,  A. Abdob, and S. X.  Ding, ``$H_-/H_{\infty}$ fault
detection filter design for discrete-time
 Takagi-Sugeno
fuzzy system,'' {\em Automatica}, vol. 49, no. 7, pp. 1996--2005,
2013.

\bibitem{newtsp2} B. S. Chen,  and C.F. Wu, ``Robust Scheduling Filter Design for a Class of
Nonlinear Stochastic Poisson Signal Systems,'' {\em IEEE Trans. Signal Process.}, vol. 63, no. 23, pp. 6245-6257, 2015.

%%%%%%%%%%%$$$$$$$
\bibitem{chench} B. S. Chen and C.   Wu,
``Robust optimal reference-tracking design method for stochastic
synthetic biology systems: T-S Fuzzy approach,'' {\em IEEE Trans.
Fuzzy Syst.}, vol. 18, no. 6, pp. 1144--1158, 2010.


%%%%%%%%$$$$
\bibitem{bschen_2008_1} B. S. Chen, Y.   Chang,  and Y.   Wang,
 ``Robust
$H_\infty$-stabilization design in gene networks under stochastic
molecular noises: Fuzzy interpolation approach,'' {\em IEEE Trans.
Fuzzy Syst.}, vol. 38, no. 1, pp. 25--42, 2008.

%%%%%%%%$$$$$
\bibitem{bschen2009} B. S. Chen, W. H. Chen, and H. L. Wu,
 ``Robust $H_2/H_\infty$ global linearization
filter design for nonlinear stochastic systems,'' {\em IEEE Trans.
Circuits Syst. I},   vol. 56, pp. 1441--1454, 2009.


\bibitem{chenchen} W. H. Chen  and B. S.  Chen,
 ``Robust stabilization design for nonlinear
stochastic system with Poisson noise via fuzzy interpolation
 method,''
{\em Fuzzy Sets  Syst.},  vol. 217, pp. 41--61, 2013.


\bibitem{newwu1} L. Chen, F. Wu, ``Almost sure exponential stability of the theta-method for stochastic differential equations,'' {\em Statistics \& Probability Letters},  vol. 82, no. 9, pp. 1669--1676, 2012.


%%%%%%%%%%$$$$$$$$$$
\bibitem{Stevenbook} S. X. Ding, \textit{Model-based Fault
Diagnosis Techniques Design Schemes, Algorithms, and Tools}.
Springer-Verlag Berlin Heidelberg, 2008.

%%%%%%%%%$$$$$
\bibitem{v.dragan_2005} V. Dragan, T. Morozan, and A. M. Stoica,
 \textit{ Mathematical
Methods in Robust Control of  Linear Stochastic Systems.} New York:
Springer, 2006.



\bibitem{newtsp3} H. Dong,  Z. Wang, S. X. Ding, and H. Gao ``Event-Based $H_{\infty}$ Filter Design for a Class of
Nonlinear Time-Varying Systems With Fading
Channels and Multiplicative Noises,'' {\em IEEE Trans. Signal Process.}, vol. 63, no. 13, pp. 3387--3395, 2015.


%%%%%%%%&&&&
\bibitem{Gershon_05}  E. Gershon,  U. Shaked,  and U. Yaesh,
 \textit{Control and Estimation of State-Multiplicative
Linear Systems}. London: Springer-Verlag, 2005.

%%%%%%%%$$$$$
\bibitem{Hinrichsen} D. Hinrichsen and A. J. Pritchard,
``Stochastic $H_{\infty}$,'' {\em SIAM J. Control Optim.}, vol. 36,
pp. 1504--1538, 1998.

\bibitem{newkao2} Y.  Kao, G. Yang, J. Xie, and L. Shi, ``$H_{\infty}$ adaptive control for uncertain Markovian jump systems with general unknown transition rates,'' {\em Applied Mathematical Modelling},  vol. 40, no. 9-10, pp. 5200--5215, 2016.
%%%%%%%%%$$$$$$
\bibitem{liusuntong} Y. Li, K.Sun, and S. Tong.
``Observer-based adaptive fuzzy fault-tolerant optimal control for
SISO  Nonlinear Systems,'' {\em IEEE  Trans. Cybern.},
 in press.


%%%%%%$$$$
\bibitem{Lih2017} H. Li, Y. Gao, P. Shi, and H. K. Lam,
``Observer-based fault detection for nonlinear systems with sensor
fault and limited communication capacity,'' {\em IEEE Trans. Autom.
Contr.}, vol. 61, no. 9, pp. 2745--2751, 2016.

\bibitem{liliu} X. Li  and   H. H. T. Liu,
 ``Characterization of
$H_-$  index  for linear time-varying systems,'' {\em Automatica},
vol. 49, no. 5, pp. 1449--1457, 2013.

\bibitem{lizhou} X. Li and K. Zhou,
 ``A time domain approach to
robust fault detection of linear time-varying systems,''
 {\em Automatica}, vol. 45, no. 1, pp. 94--102, 2009.

\bibitem{Liyan1} Y. Li, W. Zhang, and X. K. Liu,
``$H_-$ index for discrete-time stochastic systems with Markovian
jump and multiplicative noise,'' {\em Automatica}, vol. 90,  pp.
286--293, 2018.

%%%%%%$$$$$
\bibitem{yaning1}  Y. Lin, W. Zhang, and X. Jiang,
 ``Necessary and suffcient conditions for Pareto optimality of
the stochastic systems in finite horizon,'' {\em Automatica},  vol.
94, pp. 341--348, 2018.

%%%%%%%%%%$$$
\bibitem{yaning3} Y. Lin, T. Zhang, and W. Zhang,
``Pareto-based guaranteed cost control of the uncertain mean-field
stochastic systems in infinite horizon,'' {\em Automatica},
 vol. 92,
pp. 197--209, 2018.

%%%%%%%%%%%$$$$$$$$$$
\bibitem{liuwangyang} J. Liu, J. L. Wang,
and  G. H. Yang,   ``An LMI approach to minimum sensitivity analysis
with application to fault detection,'' {\em Automatica}, vol. 41,
no. 11, pp. 1995--2004, 2005.

%%%%%%%%%%$$$$

\bibitem{Deng} W. Mao, F. Deng, and A. Wan.
``Robust $H_2/H_{\infty}$ global linearization filter design for
nonlinear stochastic time-varying delay systems,'' {\em Sci.
 China-Inform. Sci.}, vol. 59. no. 3, pp. 1--17, 2016.

%%%%%%%%%%%%%%%%%$$$$$$$$$$$$
\bibitem{maobook}
X. Mao, \textit{Stochastic Differential Equations and Applications.}
2nd Edition, Horwood, 2007.

%%%%%%%%%%%%%$$$$$$$$$$$$$
\bibitem{Meskinbook} N. Meskin and K. Khorasani,
\textit{Fault Detection and Isolation Multi-Vehicle Unmanned
Systems}. Springer New York Dordrecht Heidelberg London, 2011.

%%%%%%%%%%%$$$$$
\bibitem{Petersen} I. R. Petersen, V. A. Ugrinovskii,
and A. V. Savkin, \textit{Savkin, Robust Control Design Using
$H_{\infty}$ Methods.} New York: Springer-Verlag, 2000.

%%%%%%%%%%$$$$$$$$
\bibitem{bshen_13} B. Shen, Z. Wang,  and H. Shu,
\textit{ Nonlinear Stochastic Systems with Incomplete Information:
Filtering and Control}. Springer, London, 2013.

%%%%%%%%%%$$$$$$$$$$
\bibitem{Suxiaojie2016}
X. Su, P. Shi, L. Wu, and Y. D. Song, ``Fault detection filtering
for nonlinear switched stochastic systems.'' {\em IEEE Trans. Autom.
Contr.}, vol. 61, no. 5,  pp. 1310--1315, 2016.

%%%%%%%%%%%$$$$$$$$
\bibitem{ug} V. A.  Ugrinovskii,  ``Robust $H_\infty$ control
 in the presence of stochastic
uncertainty,''  {\em Int. J. Contr.,} vol. 71, pp. 219--237, 1998.

%%%%%%%%%%%%$$$$$$$$$
\bibitem{wangzhang} H. Wang, D. Zhang, and R. Lu,
``Event-triggered $H_{\infty}$ filter design for Markovian
 jump systems
with quantization,'' {\em Nonlinear Anal-Hybri.}, vol. 28, pp.
23--41, 2018.

\bibitem{wangyang} J. L. Wang, G. H. Yang, and J. Liu,
``An LMI approach to $H_-$ index and mixed $H_-/H_\infty$   fault
detection observer design,¡± {\em Automatica}, vol. 43, no. 9, pp.
1656¨C-1665, 2007.



\bibitem{Wangzhenhua2017} Z. Wang, P. Shi, and C. C.  Lim,
``$H_-/H_{\infty}$  fault detection observer in finite
   frequency domain for linear
parameter-varying descriptor systems,'' {\em Automatica},  vol. 86,
pp. 38--45, 2017.

\bibitem{newkao1} J. Xie, Y.  Kao, and J. H. Park, ``$H_{\infty}$ performance for neutral-type Markovian switching systems with general uncertain transition rates via sliding mode control method,'' {\em Nonlinear Analysis-Hybrid Systems},  vol. 27, pp. 416--436, 2018.

\bibitem{xjxie} X. J. Xie and  N. Duan, ``Output tracking of
high-order stochastic nonlinear systems with application to
benchmark mechanical system,'' {\em IEEE Trans. Autom. Contr.,}
  vol. 55, pp.
1197--1202, 2010.

%%%%%%%%%%%$$$$$$$$
\bibitem{Yong_1999} J. Yong and  X. Y. Zhou,
\textit{ Stochastic Control: Hamiltonian Systems and HJB
Equations.}
New York: Springer,  1999.

\bibitem{newtsp1} X. Yao, L. Wu, and W. Zheng, ``Fault detection filter design for markovian jump
singular systems with intermittent measurements,'' {\em IEEE Trans. Signal Process.}, vol. 59, no. 7, pp. 3099--3109, 2011.


\bibitem{newyin1} J. Yin, S. Khoo, and Z. Man, ``Finite-time stability theorems of homogeneous stochastic nonlinear systems,'' {\em Systems \& Control Letters},  vol. 100,  pp. 6--13, 2017.

\bibitem{newyin2} J. Yin, S. Khoo, Z. Man, and X. Yu, ``Finite-time stability and instability of stochastic nonlinear systems,'' {\em Automatica},  vol. 47, no. 12, pp. 2671--2677, 2011.

\bibitem{G.Zames_81} G. Zames, ``Feedback and
optimal sensitivity: Model reference transformation, multiplicative
seminorms and approximative inverses,'' {\em IEEE Trans.  Autom.
Contr.}, vol. 26, pp. 301--320, 1981.

 \bibitem{zhangsiam}
W. Zhang, and B. S. Chen,
 `` State feedback $H_{\infty}$  control for
a class of nonlinear stochastic systems," {\em SIAM J. Control
Optim.}, vol. 44, no. 6, pp. 1973--1991, 2006.
%
\bibitem{zhangfilter}
W. Zhang, and B. S. Chen, ``Robust $H_{\infty}$  filtering for
nonlinear stochastic systems," {\em IEEE Trans. Autom. Contr.}, vol.
53, no. 2, pp. 589--598, 2005.

\bibitem{Zhangremark2014}
 W. Zhang, B. S. Chen, H. Tang, et al,
 ``Some remarks on general nonlinear stochastic
 $H_{\infty}$ control with state, control,
 and disturbance-dependent noise,''
  {\em IEEE Trans. Autom. Contr.},
  vol. 59, no. 1, pp. 237--242, 2014.


%%%%%%%%%%%%%$$$$$$$$$
\bibitem{book}
W. Zhang, L. Xie, and B. S. Chen,     \textit{ Stochastic
$H_2/H_{\infty}$ Control: A Nash Game Approach.}  Boca Raton, FL,
USA: CRC Press, 2017.




\bibitem{zhzhang} Z. H. Zhang and G. H. Yang,
``Interval observer-based fault isolation for  discrete-time fuzzy
interconnected systems  with unknown interconnections,'' {\em IEEE
Trans. Cybern.}, vol. 47, no. 9, pp. 1--12,  2017.

\bibitem{Zhang2006}
 W. Zhang, H. Zhang, and B. S. Chen,
  ``Stochastic $H_2/H_{\infty}$ control
   with $(x,u,v)$-dependent noise:
   Finite horizon case,'' {\em Automatica},
   vol. 42, no. 11, pp. 1891--1898, 2006.

\bibitem{Zhong2003}  M. Zhong, S. X. Ding,
J. Lam, and H. Wang, ``An LMI approach to design robust fault
detection filter for uncertain LTI systems,'' {\em Automatica}, vol.
39, no. 3, pp. 543--550, 2003.







\end{thebibliography}
\end{document}